\documentclass[final,leqno,onefignum,onetabnum,dvips,
letterpaper=true,colorlinks=true,linkcolor=red,
filecolor=green,citecolor=red,pdfpagemode=None]{siamltex}

\usepackage{hyperref}
\usepackage{array,arydshln}
\usepackage{times}
\usepackage{wrapfig}
\usepackage{subcaption}
\usepackage{graphicx}
\usepackage{amsmath}
\usepackage{amssymb,thmtools}
\usepackage{cancel}
\usepackage{floatrow}
\usepackage{enumitem}
\usepackage{chngcntr}
\usepackage{multirow}
\usepackage{pdflscape}
\usepackage{setspace}
\usepackage{blkarray}
\usepackage{relsize}
\usepackage{xcolor}
\usepackage{framed}
\usepackage{mdframed}
\usepackage{stackrel}
\usepackage[font=small,labelfont=bf]{caption}
\usepackage[titletoc,title]{appendix}
\usepackage{soul}
\usepackage[normalem]{ulem}

\newcommand{\bflambda}  { {\boldsymbol \lambda} }

\newcommand{\bftheta}   { {\boldsymbol \theta} }

\newcommand{\bfLambda}  { {\boldsymbol \Lambda} }
\newcommand{\bfTheta}   { {\boldsymbol \Theta} }

\newcommand{\zero}  { {\bf 0} }

\newcommand{\hs}{\hskip 0.1em }
\newcommand{\partialsp}{\partial\hskip 0.1em }

\newcommand{\bfd}{{\bf d}}
\newcommand{\bfe}{{\bf e}}
\newcommand{\bff}{{\bf f}}
\newcommand{\bfg}{{\bf g}}

\newcommand{\bfk}{{\bf k}}

\newcommand{\bfm}{{\bf m}}
\newcommand{\bfn}{{\bf n}}

\newcommand{\bfq}{{\bf q}}
\newcommand{\bfr}{{\bf r}}

\newcommand{\bft}{{\bf t}}
\newcommand{\bfu}{{\bf u}}
\newcommand{\bfv}{{\bf v}}
\newcommand{\bfw}{{\bf w}}
\newcommand{\bfx}{{\bf x}}
\newcommand{\bfy}{{\bf y}}
\newcommand{\bfz}{{\bf z}}

\newcommand{\bfA}{{\bf A}}
\newcommand{\bfB}{{\bf B}}
\newcommand{\bfC}{{\bf C}}
\newcommand{\bfD}{{\bf D}}

\newcommand{\bfF}{{\bf F}}

\newcommand{\bfH}{{\bf H}}
\newcommand{\bfI}{{\bf I}}
\newcommand{\bfJ}{{\bf J}}

\newcommand{\bfL}{{\bf L}}

\newcommand{\bfN}{{\bf N}}

\newcommand{\bfQ}{{\bf Q}}

\newcommand{\bfT}{{\bf T}}

\newcommand{\bfV}{{\bf V}}
\newcommand{\bfW}{{\bf W}}

\newcommand{\bfY}{{\bf Y}}

\newcommand{\opM}{{\mathrm{M}}}

\newcommand{\opR}{{\mathrm{R}}}

\newcommand*\fd[1]{{#1}_\tau}

\newcommand*{\fdn}{{\fd{\bfn}}}

\newcommand*{\fdy}{{\fd{\bfy}}}

\newcommand*{\qed}{\hfill\ensuremath{\square}}

\newcommand{\mdiag}[1]{\mathrm{diag}\left( #1 \right)}

\counterwithin{table}{section}

\setlist[itemize]{leftmargin=1em}
\setlist[enumerate]{leftmargin=1.5em}

\newcolumntype{L}[1]{>{\raggedright\let\newline\\\arraybackslash\hspace{0pt}}m{#1}}
\newcolumntype{C}[1]{>{\centering\let\newline\\\arraybackslash\hspace{0pt}}m{#1}}
\newcolumntype{R}[1]{>{\raggedleft\let\newline\\\arraybackslash\hspace{0pt}}m{#1}}

\newtheorem{algorithm}[theorem]{Algorithm}

\definecolor{turquoise}{cmyk}{0.65,0,0.1,0.1}
\definecolor{purple}{rgb}{0.65,0,0.65}
\definecolor{green}{rgb}{0, 0.5, 0}
\definecolor{blue}{rgb}{0, 0, 1}
\definecolor{orange}{rgb}{0.8, 0.6, 0.2}
\definecolor{red}{rgb}{0.8, 0.2, 0.2}
\definecolor{brown}{rgb}{0.5, 0.16, 0.16}

\title{The Discrete Adjoint Method for Exponential Integration}

\author{Kai Rothauge\thanks{Department of Mathematics, University of British Columbia (\tt{rothauge@math.ubc.ca})} \and Eldad Haber\thanks{Department of Earth and Ocean Science, University of British Columbia (\tt{haber@eos.ubc.ca})} \and Uri Ascher\thanks{Department of Computer Science, University of British Columbia (\tt{ascher@cs.ubc.ca})} }

\begin{document}
\maketitle


\begin{abstract}
The {implementation} of the discrete adjoint method {for} exponential time differencing (ETD) schemes 
is {considered}. This is important for parameter estimation problems that are constrained by 
stiff time-dependent PDEs when the discretized PDE system is solved using an exponential integrator. 
We also discuss the closely related topic of computing the action of the sensitivity matrix on a vector, 
which is required when performing a sensitivity analysis. The PDE system is assumed to be semi-linear 
and can be the result of a linearization of a nonlinear PDE, leading to exponential Rosenbrock-type methods. 
We discuss the computation of the derivatives of the $\varphi$-functions that are used by ETD schemes 
and find that the derivatives strongly depend on the way the $\varphi$-functions are evaluated numerically. 
A general adjoint exponential integration method, required when computing the gradients, is developed 
and its implementation is illustrated by applying it to the Krogstad scheme. 
The applicability of the methods developed here to pattern formation problems is demonstrated 
using the Swift-Hohenberg model.
\end{abstract}



\pagestyle{myheadings}
\thispagestyle{plain}



\section{Introduction}
\label{sec:introduction}

Large-scale {distributed} parameter estimation is an important type of inverse problem dealing with the recovery or approximation of the model parameters appearing in partial differential equations (PDEs). 
The PDE system models some real-life process and observations of this process, which usually contain noise, 
are compared to numerically computed solutions of the PDEs. 
Parameter estimation problems, often referred to as model calibration, are very common in engineering, many branches of science, economics, and elsewhere.

In this paper we consider such inverse problems in the context of stiff, time-dependent, semi-linear PDEs. Using the method of lines approach, the PDE is discretized in space by applying a finite difference, finite volume or finite element method, where the (known) boundary and initial conditions are assumed to have already been incorporated. This results in a large system of ordinary differential equations (ODEs) that can be written in generic first-order form as
\begin{equation}
\begin{aligned}
	\dfrac{\partialsp\bfy(t;\bfm)}{\partialsp t} &= \bff(\bfy(t;\bfm),t,\bfm) \\
&= \bfL(\bfm)\hs\bfy(t;\bfm) + \bfn(\bfy(t;\bfm),t,\bfm) \\
	\bfy(0) &= \bfy_0,
\end{aligned}\label{eqn:general_semidiscrete_nonlinear_system}
\end{equation}
on an underlying discrete spatial grid, with $0 \leq t \leq T$. The vector $\bfm \in \mathbb{R}^{N_{\bfm}}$ is the set of discretized model parameters that could, for instance, represent some physical properties of the underlying material that we want to estimate. We refer to $\bfy = \bfy(t;\bfm)$ in this context as the \textit{forward solution}. It is a time-dependent vector of length $N$ and depends on $\bfm$ indirectly through the discretized operators $\bfL$ and $\bfn$. The operator $\bfL$ contains the leading derivatives and is linear, while $\bfn$ is generally nonlinear in $\bfy$. Many interesting PDEs have this form.

The problem of estimating the parameters in \eqref{eqn:general_semidiscrete_nonlinear_system}
based on a given set of measurements is often tackled using gradient-based optimization procedures. For large-scale problems, the adjoint method is a well-known efficient approach to computing the required gradients. The goal of this study is to systematically implement the discrete adjoint method for the case where a fully discretized version of \eqref{eqn:general_semidiscrete_nonlinear_system}
is solved using an exponential time differencing (ETD) 
scheme \cite{HochbruckOstermann2010, CoxMatthews2002, Krogstad2005}, 
also known as an exponential integrator.
These methods have come to form an important approach for numerically solving PDEs,
particularly when the PDE system exhibits stiffness in time;
yet they pose several challenges in carrying out their corresponding adjoint method.  
Our results are also relevant to the problem of sensitivity analysis, where the efficient computation of the action of the sensitivity matrix on a vector is required.

\subsection{PDE-Constrained Optimization}

The parameter estimation problem is to recover a feasible set of parameters $\bfm$ that approximately solves the 
PDE-constrained optimization problem
\begin{equation}
	\bfm^\star = \underset{\bfm}{\arg\min}\hs\hs \Omega(\bfm)  \qquad \mathrm{s.t. \hs\hs \eqref{eqn:general_semidiscrete_nonlinear_system} \hs\hs holds},
	\label{eqn:minimize_Omega}
\end{equation}
where the \textit{objective function}
$$
\Omega(\bfm) = \Omega(\bfd, \bfd^\mathrm{obs}; \bfm, \bfm^\mathrm{ref}) = 
\opM(\bfd, \bfd^\mathrm{obs}) + \beta\hs\opR(\bfm, \bfm^\mathrm{ref})
$$
has two components. The \textit{regularization function} $\opR(\bfm, \bfm^\mathrm{ref})$ 
penalizes straying too far away from the prior knowledge we have of the true parameter values 
(which is incorporated in $\bfm^\mathrm{ref}$). Common examples are Tikhonov-type regularization~\cite{EnglHankeNeubauer1996,vandenDoelAscherHaber2013},
including least squares and total variation (TV)~\cite{Vogel2002}. 
The \textit{data misfit function} $\opM(\bfd, \bfd^\mathrm{obs})$ in some way quantifies 
the difference between $\bfd$ and $\bfd^\mathrm{obs}$, 
with $\bfd^\mathrm{obs}$ a given set of observations of the true solution and 
$\bfd = \bfd(\bfm) = \bfd(\bfy(\bfm))$ the observation of the current simulated solution, 
depending on $\bfm$ implicitly through the forward solution. 

The most popular choice for $\opM$, corresponding to the assumption that the noise 
in the data is simple and white, 
is the least-squares function $\opM = \frac{1}{2} \left\|\hs\bfd - \bfd^\mathrm{obs}\hs\right\|_2^2$; 
but we will not restrict our discussion to any particular misfit function. 
The relative importance of the data misfit and prior 
is adjusted using the regularization parameter $\beta$. 

Several classes of optimization procedures 
can be used to solve~\eqref{eqn:minimize_Omega}, 
see for instance \cite{EnglHankeNeubauer1996,Vogel2002,vandenDoelAscherHaber2013}.
Typically, all reduced space methods for PDE-constrained optimization
require the gradient of $\Omega$ with respect to $\bfm$, 
$\nabla_{\bfm}\Omega=\nabla_{\bfm}\opM + \beta\hs\nabla_{\bfm}\opR$ \cite{NocedalWright2006,DennisSchnabel1996}. 
The derivatives of the regularization function are known and are independent of the time-stepping scheme, 
so in this paper we focus exclusively on computing the gradient of the misfit function $\opM$.

\subsection{The Adjoint Method}

The gradient of the misfit function requires computing the action of $\bfJ^T$, where
$\bfJ = \frac{\partial\bfd}{\partial\bfm}$ is the sensitivity matrix. 
This matrix function  stores the first-order derivatives of the predicted data with respect to the model parameters and can be calculated explicitly in small-scale applications. However, in cases such as those considered here it is far more feasible to compute the action of the 
sensitivity matrix on a vector using the adjoint method. Originally developed in the optimal control community, 
the adjoint method was introduced in~\cite{Chavent1974} to the theory of inverse problems 
to efficiently compute the gradient of a function. See \cite{Plessix2006} for a review of the method applied to geophysical problems.

There are two frameworks that one can take when computing derivatives of the misfit function, discretize-then-optimize (DO) or optimize-then-discretize (OD). In OD one forms the adjoint differential problem, for either the given PDE or its semi-discretized form~\eqref{eqn:general_semidiscrete_nonlinear_system}, 
which is subsequently discretized, thus decoupling the discretization process of the problem from that of its adjoint. A disadvantage of this approach is that the gradients of the discretized problem are not obtained exactly, because a discretization error gets introduced when moving from the continuous setting~\cite{Gunzburger2002}. 
We therefore prefer the DO approach, where one applies the adjoint method to the fully discretized form \eqref{eqn:general_semidiscrete_nonlinear_system}. 
The time-stepping method is therefore of central importance and this article addresses the implementation of the discrete adjoint method where the time-stepping method is an ETD scheme.

An alternative approach to computing the gradient that falls inside the DO framework 
is automatic differentiation (AD) \cite{GriewankWalther2008,Neidinger2010}.
AD is a great tool for many purposes, but for large-scale problems where efficiency and memory allocation are essential, the adjoint method approach is advantageous when performance is crucial. In such circumstances it is preferable to have hard-coded and optimizable gradient and sensitivity computations. Having an explicit algorithm for the adjoint time-stepping method, as we develop here, allows one to exploit opportunities to increase the computational efficiency, such as parallelization and the precomputation of repeated quantities. Furthermore, it would require a very sophisticated AD code to compute the derivatives of the 
$\varphi$-functions that arise in ETD schemes, further elaborated upon in Sections \ref{sec:ETD} and \ref{sec:action_of_phi}.

\subsection{Main Contributions and Structure of this Article}

Exponential integration methods are reviewed in Section \ref{sec:ETD}. 
It is shown how to abstractly represent them in the form $\bft(\bfy,\bfm) = \zero$, which we need when applying the discrete adjoint method. We also review Rosenbrock-type methods, since our techniques can be applied to them as well 
\cite{HochbruckVanDenEshof2006a, HochbruckOstermannSchweitzer2009, LuanOstermann2014a, LuanOstermann2016}. 

ETD schemes require the use of $\varphi$-functions that are introduced in Section~\ref{sec:ETD} 
and whose numerical evaluation is briefly reviewed in Section~\ref{sec:action_of_phi}.

In Section \ref{sec:discrete_adjoint_method} we apply the discrete adjoint method to the 
abstract time-stepping representation derived in Section \ref{sec:ETD}. Derivatives of $\bft$ with respect to both the solution and the model parameters $\bfm$ are required, and we derive expressions for them in Section \ref{sec:derivatives_of_T}. 
These expressions require
the derivatives of the 
$\varphi$-functions when the discretized linear operator $\bfL$ depends on the solution or the model parameters, and we discuss their differentiation in Section \ref{sec:derivatives_of_phi}.

Section \ref{sec:linearized_forward_problem} discusses the linearized forward problem that 
must be solved when computing the action of the sensitivity matrix on some vector.
In Section \ref{sec:adjoint_problem} we present the solution of the adjoint problem needed 
when computing the action of the transpose of the sensitivity matrix, 
thereby also playing a central role in the calculation of the gradient $\nabla_\bfm \opM$. 
We give general algorithms applicable to any ETD scheme in these sections, 
and also illustrate their implementation using the ETD scheme of~\cite{Krogstad2005}.

The results of this work are used in Section \ref{sec:numerical_example}
 to solve a simple parameter estimation problem involving the Swift-Hohenberg model, a PDE problem with solutions exhibiting the 
interesting phenomenon of pattern formation. 
Some final thoughts and avenues for future work are provided in Section \ref{sec:conclusion}. 



\section{Exponential Time Differencing}
\label{sec:ETD}

When considering time discretizations for the semi-discretized PDE~\eqref{eqn:general_semidiscrete_nonlinear_system}, 
it is often convenient to suppress the dependence on $\bfm$ in the notation.

Exponential time differencing methods, also referred to as exponential integrators, 
were originally developed in the 1960s \cite{Certaine1960,Pope1963}
and have attracted much recent attention~\cite{HochbruckOstermann2010, LuanOstermann2014b, LuanOstermann2013, LuanOstermann2014c, Tokman2011, Tokman2005}.
This has become an important approach to numerically solving PDEs, particularly when the PDE system exhibits stiffness. 
Let us discretize the time interval 
as $0 = t_0 < t_1 < \cdots < t_K = T$, 
with a possibly variable time-step size $\tau_k = t_{k+1} - t_{k}$. Consider \eqref{eqn:general_semidiscrete_nonlinear_system} and let $\bfy_k$ denote the approximation
of its solution $\bfy(t_k)$ at time $t_k$.

Before we review exponential integration, we first give a quick overview of Rosenbrock-type methods. Here, an ODE system resulting from a spatial semi-discretization of 
a fully nonlinear time-dependent PDE $\frac{\partialsp}{\partialsp t}\hs\bfy(t) = \bff(\bfy(t),t)$  
is written in semi-linear form, after which an exponential integrator can be applied. This can be done by writing
$$
	\bff(\bfy(t),t) = \bfL_k\hs\bfy(t) + \bfn_k(\bfy(t),t),
$$
where $\bfn_k\left(\bfy(t),t\right) = \bff(\bfy(t),t) - \bfL_k\hs\bfy(t)$.
ETD schemes are applied to this linearization, 
involving exponentiation of the matrix $\bfL_k$ and
leading to exponential Rosenbrock-type methods \cite{HochbruckVanDenEshof2006a, HochbruckOstermannSchweitzer2009}.
The matrix $\bfL_k$ is meant to approximate the Jacobian $\dfrac{\partialsp\bff}{\partialsp\bfy}(\bfy_k,t_k)$
at the $k$th time-step.
Occasionally it is sufficient to select $\bfL_k=\bfL$ independently of time, but when the dynamical system trajectory
varies significantly and rapidly we may have to set $\bfL_k = \bfL_k\left(\bfy_k\right) = \dfrac{\partialsp\bff}{\partialsp\bfy}(\bfy_k,t_k)$.
The case where the linear operator $\bfL_k$ depends on $\bfy_k$ 
adds a significant amount of complexity to the adjoint computations considered later in this paper.


Exponential integration is briefly reviewed in Section \ref{subsec:derivation_of_ETDRK} and we give examples of particular schemes in Section \ref{subsec:ETDRK_examples}. The discrete adjoint method discussed in Section \ref{sec:discrete_adjoint_method} requires that the time-stepping method be abstractly represented in the form of a \textit{discrete time-stepping equation} 
\begin{equation}
	\bft\left(\fdy,\bfm\right) = \zero,	\label{eqn:abstract_timestepping_system}
\end{equation}
where $\bft$ is a vector representing the time-stepping method and $\fdy$ denotes
the approximate solution vector composed of all the $\bfy_k$'s.
In Section \ref{subsec:representation_by_T} 
we show how to represent an arbitrary ETD scheme in the form \eqref{eqn:abstract_timestepping_system}.



\subsection{Derivation of ETDRK methods}
\label{subsec:derivation_of_ETDRK} 

Recall the semi-discretized, semi-linear system \eqref{eqn:general_semidiscrete_nonlinear_system} 
which we write as
\begin{equation}
 \dfrac{\partialsp\bfy(t)}{\partialsp t} = \bfL_k\hs\bfy(t) + \bfn_k\left(\bfy(t),t\right).
\label{eqn:general_semidiscrete_semilinear_system}
\end{equation} 
Integrating \eqref{eqn:general_semidiscrete_semilinear_system} exactly from time level $t_k$ to $t_{k+1} = t_k + \tau_{k}$ 
gives
\begin{equation}
\bfy_{k+1} = e^{\tau_{k}\bfL_k}\bfy_k + \int_{t_k}^{t_{k+1}}e^{(t_{k+1}-t)\bfL_k}\,\bfn_k\left(\bfy(t_k+t),t_k+t\right)dt.
\label{eqn:Volterra_integral_equation}
\end{equation}
The \textit{exponential Euler method} is obtained by  interpolating the 
integrand at the known value $\bfn_k\left(\bfy_k,t_k\right)$ only, 
\begin{equation}
\bfy_{k+1} = e^{\tau_{k}\bfL_k}\bfy_k + \tau_{k}\varphi_1(\tau_{k}\bfL_k)\bfn_k\left(\bfy_k,t_k\right)
\label{eqn:exponential_Euler_method}
\end{equation}
where $\varphi_1(z) = \frac{e^z - 1}{z}$. This is the simplest numerical method that can be obtained for solving \eqref{eqn:Volterra_integral_equation}.

The integral in  \eqref{eqn:Volterra_integral_equation} can be approximated using some quadrature rule, leading to a class of $s$-stage explicit \textit{exponential time differencing Runge-Kutta} (ETDRK) methods with matrix coefficients $a_{i,j}(\tau_k{\bfL}_k)$, weights $b_{i}(\tau_k{\bfL}_k)$ and nodes $c_{i}$, so for $1\leq i,j\leq s$ we obtain
\begin{subequations}
\label{eqn:ETDRK_general_method}
\begin{equation}
	\bfy_{k+1} = e^{\tau_k{\bfL}_k}\bfy_k + \tau_{k}\sum^s_{i=1}\,b_i (\tau_k{\bfL}_k)\bfY_{k+1,i},
\label{eqn:ETDRK_general_formula}
\end{equation}
with the internal stages 
\begin{equation}
	\bfY_{k+1,i} = 
	\bfn_k\left(e^{c_i\tau_k{\bfL}_k}\bfy_{k} + \tau_{k}\sum^{i-1}_{j=1}a_{ij}(\tau_k{\bfL}_k)\bfY_{k+1,j},t_k + c_i\tau_{k}\right) \qquad 1\leq i\leq s.
\label{eqn:ETDRK_general_internal_stages}
\end{equation}
\end{subequations}
The procedure starts from a known initial condition $\bfy_0$. There are several alternate ways of writing \eqref{eqn:ETDRK_general_formula}, but for our purposes the representation given here is most useful.

The Butcher tableau for these methods is
$$
\begin{array}{c|ccccc}
c_1    &     \\
c_2    & a_{21}(\tau{\bfL}_k) &     \\
\vdots & \vdots & \ddots   \\
c_s    & a_{s1}(\tau_k{\bfL}_k) & \cdots & a_{s,s-1}(\tau_k{\bfL}_k) &   \\
\hline \\[-1em]
       & b_1(\tau_k{\bfL}_k) & \cdots & b_{s-1}(\tau_k{\bfL}_k) & b_s(\tau_k{\bfL}_k)
\end{array}.
$$
The coefficients $a_{ij}$ and $b_{i}$ are linear combinations of the entire functions
$$
	\varphi_0(z) = e^z \qquad \varphi_\ell = \int_0^1 e^{(1-\theta)z}\frac{\theta^{\ell-1}}{(\ell-1)!}\,d\theta, \quad \ell \geq 1.
$$
It is not hard to see that the $\varphi$-functions satisfy the recurrence relation
\begin{equation}
	\varphi_{\ell}(z) = \dfrac{\varphi_{\ell-1}(z) - \varphi_{\ell-1}(0)}{z}, \quad \ell > 0,
	\label{eqn:phi_recurrence_relation}
\end{equation}
and that $\varphi_{\ell}(z) = \sum_{i=0}^{\infty}\dfrac{z^i}{(i+\ell)!}$. Notice that the expansion of $\varphi_{\ell}(z)$ is that of the exponential function with the coefficients shifted forward.

As is evident by the structure of $\varphi_\ell$ for $\ell > 0$, for small $z$ the evaluation of $\varphi_\ell (z)$ will be subject to cancellation error, and this could become a problem when evaluating $\varphi_\ell (\tau_k{\bfL}_k)$ if the matrix $\tau_k{\bfL}_k$ has small eigenvalues.


From now on, for brevity of notation we use $\varphi_{\ell} = \varphi_{\ell}(\tau_k{\bfL}_k)$ and $\varphi_{\ell,i} = \varphi_{\ell}(c_{i}\tau_k{\bfL}_k)$. 


\subsection{Examples}
\label{subsec:ETDRK_examples} 

The four-stage ETD4RK method of Cox and Matthews \cite{CoxMatthews2002} has the following Butcher tableau:
\begin{equation}
\begin{array}{c|ccccc}
0    &   \\
\frac{1}{2} & \frac{1}{2}\varphi_{1,2}  \\
\frac{1}{2} & \zero_{N\times N} & \frac{1}{2}\varphi_{1,3} &     \\
1    & \varphi_{1,4} - \varphi_{1,3} & \zero_{N\times N} & \varphi_{1,3} &   \\
\hline \\[-1em]
       & \varphi_{1} - 3\varphi_{2} + 4\varphi_{3} & 2\varphi_{2} - 4\varphi_{3} & 2\varphi_{2} - 4\varphi_{3} & 4\varphi_{3} - \varphi_{2}
\end{array}.
\label{eqn:ETDRK_Butcher_tableau_Cox_Matthews}
\end{equation}
This method can be fourth-order accurate when certain conditions are satisfied, but in the worst case is only second-order.

Krogstad \cite{Krogstad2005} derived the method given by
\begin{equation}
\begin{array}{c|ccccc}
0    &     \\
\frac{1}{2} & \frac{1}{2}\varphi_{1,2} \\
\frac{1}{2} & \frac{1}{2}\varphi_{1,3} - \varphi_{2,3} & \varphi_{2,3} &     \\
1    & \varphi_{1,4} - 2\varphi_{2,4} & \zero_{N\times N} & 2\varphi_{2,4} &   \\
\hline \\[-1em]
       & \varphi_{1} - 3\varphi_{2} + 4\varphi_{3} & 2\varphi_{2} - 4\varphi_{3} & 2\varphi_{2} - 4\varphi_{3} & 4\varphi_{3} - \varphi_{2}
\end{array}.
\label{eqn:ETDRK_Butcher_tableau_Krogstad}
\end{equation}
It is usually also fourth-order accurate and has order three in the worst case.

The following five-stage method is due to Hochbruck and Ostermann~\cite{HochbruckOstermann2005a}:
\begin{equation}
\begin{array}{c|ccccc}
0    &     \\
\frac{1}{2} & \frac{1}{2}\varphi_{1,2} \\
\frac{1}{2} & \frac{1}{2}\varphi_{1,3} - \varphi_{2,3} & \varphi_{2,3} \\
1    & \varphi_{1,4} - 2\varphi_{2,4} & \varphi_{2,4} & \varphi_{2,4}  \\
\frac{1}{2} & \frac{1}{2}\varphi_{1,5}- \frac{1}{4}\varphi_{2,5} - a_{5,2} & a_{5,2} & a_{5,2} & \frac{1}{4}\varphi_{2,5} - a_{5,2} &   \\
\hline \\[-1em]
       & \varphi_{1} - 3\varphi_{2} + 4\varphi_{3} & \zero_{N\times N} & \zero_{N\times N} & -\varphi_{2} + 4\varphi_{3} & 4\varphi_{2} - 8\varphi_{3}
\end{array},
\label{eqn:ETDRK_Butcher_tableau_Hochbruck_Ostermann}
\end{equation}
with $a_{5,2} = \frac{1}{2}\varphi_{2,5} - \varphi_{3,4} + \frac{1}{4}\varphi_{2,4} - \frac{1}{2}\varphi_{3,5}$. It has order four under certain mild assumptions.


\subsection{Representation by $\bft$}
\label{subsec:representation_by_T} 

The method given by \eqref{eqn:ETDRK_general_method} can be abstractly represented by
\begin{equation}
    \begin{bmatrix}
         \bfI_N  &       \\
                 & \bfI_{sN}  \\
        -e^{\tau_k{\bfL}_k} & -\bfB_{k+1}^\top & \bfI_N 
    \end{bmatrix}
    \begin{bmatrix}
        \bfy_k     \\
        \bfY_{k+1}  \\
        \bfy_{k+1}
    \end{bmatrix} - \begin{bmatrix}
        \zero_{N \times 1} \\
        \bfN_{k+1}     \\
        \zero_{N \times 1}
    \end{bmatrix} = \begin{bmatrix}
        \bfy_k     \\ 
        \zero_{sN \times 1} \\
        \zero_{N \times 1}
    \end{bmatrix},
    \label{eqn:ETDRK_single_step_matrix}
\end{equation}
with
\begin{subequations}
\begin{equation}
    \bfY_{k+1} = \begin{bmatrix}
         \bfY_{k+1,1} \\
         \vdots \\
         \bfY_{k+1,s}
    \end{bmatrix}, \quad \bfB_{k+1} = \tau_{k}\begin{bmatrix}
         b_1(\tau_k{\bfL}_k)^\top \\
         \vdots \\
         b_s(\tau_k{\bfL}_k)^\top
    \end{bmatrix} \quad \mathrm{and} \quad \bfN_{k+1} = \begin{bmatrix}
         \bfn_{k,1} \\
         \vdots \\
         \bfn_{k,s}
    \end{bmatrix} ,
	\label{eqn:ETDRK_assorted_stuff}
\end{equation}
where
\begin{equation}
\bfn_{k,i} = \bfn_k\left(e^{c_i\tau_k{\bfL}_{k}}\bfy_{k} + \tau_{k}\sum\limits^{i-1}_{j=1}a_{ij}(\tau_k{\bfL}_{k})\bfY_{k+1,j},t_{k} + c_{i}\tau_{k}\right).
\end{equation}
\end{subequations}

In~\eqref{eqn:ETDRK_single_step_matrix} we have a single time step in the solution procedure, so the ETDRK procedure as a whole can be represented by
\begin{equation}
	\bft\left(\fdy\right) = \bfT\hs\fdy - \fdn\left(\fdy\right) - \bfq = \zero,
\label{eqn:ETDRK_T}
\end{equation}
with
\begin{equation}
\begin{aligned}
\fdy &= \begin{bmatrix} \bfY_1^\top & \bfy_1^\top & \bfY_2^\top & \bfy_2^\top & \cdots & \bfY_K^\top & \bfy_K^\top \end{bmatrix}^\top , \\
\fdn &= \begin{bmatrix} \bfN_1^\top & \zero_{1\times N} & \bfN_2^\top & \zero_{1\times N} & \cdots & \bfN_{K}^\top & \zero_{1\times N} \end{bmatrix}^\top ,\\
\bfq &= \begin{bmatrix} \zero_{1\times sN} & (e^{\tau_0{\bfL}_0}\bfy_0)^\top & \zero_{1\times (s+1)N} & \cdots & \zero_{1\times (s+1)N} \end{bmatrix}^\top , \\
	\bfT &= \begin{bmatrix}
   \bfI_{sN} \\
   -\bfB_{1}^\top  &  \bfI_N  \\
                &         & \bfI_{sN} \\
                & -e^{\tau_1{\bfL}_1} & -\bfB_{2}^\top & \bfI_N  \\
                &         &  \ddots      & \ddots & \ddots  \\
                &         &              &       &  & \bfI_{sN} \\
                &         &              &        &   -e^{\tau_{K-1}{\bfL}_{K-1}} & -\bfB_{K}^\top & \bfI_N
\end{bmatrix},
\end{aligned}
\label{eqn:ETDRK_miscellaneous_definitions}
\end{equation}
which is a block lower-triangular matrix. We have explicitly included the internal stages in $\fdy$ because they will be needed later on.



\section{The Action of $\varphi_\ell$ on Arbitrary Vectors}
\label{sec:action_of_phi}

To simplify the notation in this section, let $\bfL = \tau_k{\bfL}_k$ and $\varphi = \varphi_\ell$ for some $k,\ell \geq 0$. We very briefly review some methods used in practice to evaluate the product of $\varphi(\bfL) \in \mathbb{R}^{N\times N}$ with some arbitrary vector $\bfw \in \mathbb{R}^N$, where $N$ is too large for $\varphi(\bfL)$ to be computed explicitly and stored in full, or where it is impractical to first diagonalize $\bfL$. 

Much has been written about the approximation of these products for large $N$;
see~\cite{HochbruckOstermann2010} and the references therein. Here we mention four of the most relevant approaches to performing these approximations, all of which also help to address the numerical cancellation error that would occur when computing $\varphi$ directly using the recurrence relation \eqref{eqn:phi_recurrence_relation}.

If using a Rosenbrock-type scheme with $\bfL$ depending on $\bfy_k$, 
then we are interested in finding methods that lend themselves to 
calculating the derivatives of $\varphi(\bfL(\bfy_k,\bfm))\bfw$ with respect to $\bfy_k$ or $\bfm$. This will be explored in more detail in Section \ref{sec:derivatives_of_phi}.


\subsection{Krylov Subspace Methods}
\label{subsec:Krylov_subspace_methods}

The $M$th Krylov subspace with respect to a matrix $\bfL$ and a vector $\bfw$ is denoted by $\mathcal{K}_M(\bfL,\bfw) = \mathrm{span}\{\bfw,\bfL\bfw,\ldots,\bfL^{M-1}\bfw\}$.
Normalizing $\left\|\bfw\right\|=1$, the Arnoldi process can be used to construct an orthonormal basis $\bfV_M \in \mathbb{C}^{N \times M}$ of $\mathcal{K}_M(\bfL,\bfw)$ and an unreduced upper Hessenberg matrix $\bfH_M \in \mathbb{C}^{M \times M}$ satisfying the standard Krylov recurrence formula
$$
	\bfL\,\bfV_M = \bfV_M\bfH_M + h_{M+1,M}\bfv_{M+1}\bfe^T_M,\qquad \bfV_M^*\bfV_M = \bfI_M,
$$
with $\bfe_M$ the $M$th unit vector in $\mathbb{C}^N$. 
Using the orthogonality of $\bfV_M$, it can then be shown that
\begin{equation}
	\varphi(\bfL)\bfw \approx \bfV_M\varphi(\bfH_M)\bfe_1 
	\label{eqn:phi_product_Krylov}
\end{equation}
(see for instance~\cite{Saad2003}).
It is assumed that $M \ll N$, so that $\varphi(\bfH_M)$ can be computed using standard methods such as diagonalization or Pad{\'e} approximations.

There has been a lot of work on Krylov subspace methods for evaluating matrix functions, see for instance
\cite{EiermannErnst2006,Guttel2013} and the references therein.
See in particular~\cite{HochbruckLubich1997} for a discussion on Krylov subspace methods for matrix exponentials.


\subsection{Polynomial Approximations}
\label{subsec:polynomial_approximation} 

Polynomial methods approximate $\varphi(\bfL)$ using some truncated polynomial series, for instance Taylor series (which is rarely used in this context), Chebyshev series for Hermitian or skew-Hermitian $\bfL$, Faber series for general $\bfL$, or Leja interpolants. See \cite{HochbruckOstermann2010} and the references therein for a review of Chebyshev approximations and Leja interpolants in the context of exponential time differencing.

A polynomial approximation can generally be written in the form 
\begin{equation}
	\varphi(\bfL)\bfw \approx \sum_{j=0}^M\,c_j\,\bfL^j\bfw,
	\label{eqn:phi_product_general_polynomial}
\end{equation}
although sometimes other forms are more suitable. 
For instance, in the case of Chebyshev polynomials it makes more sense to write
\begin{equation}
	\varphi(\bfL)\bfw \approx \sum_{j=0}^M\,c_j\,T_j(\bfL)\bfw
	\label{eqn:phi_product_chebyshev}
\end{equation}
if $\bfL$ is Hermitian or skew-Hermitian and the eigenvalues of $\bfL$ all lie inside $[-1,1]$. The $T_j(\bfL)$ satisfy the recurrence relation
$$
	T_{j+1}(\bfL) = 2\bfL\,T_{j}(\bfL) - T_{j-1}(\bfL), \qquad j = 1,2,\ldots
$$
initialized by $T_{0}(\bfL) = \bfI$ and $T_{1}(\bfL) = \bfL$. 


\subsection{Rational Approximations}
\label{subsec:rational_approximation} 

The function $\varphi(z)$ can be estimated to arbitrary order using rational approximations
$$
	\varphi(z) \approx \varphi_{\left[m,n\right]}(z) = \sum_{i=0}^{m}a_iz^i\,/\,\sum_{k=0}^{n}b_kz^k = p_m(z)/q_n(z).
$$
The polynomials $p_m(z)$ and $q_n(z)$ can be found using either Pad{\'e} approximations or by using the Carath{\'e}odory-Fej{\'e}r (CF) method on the negative real line, which is an efficient method for constructing near-best rational approximations. 
It has been applied to the problem of approximating $\varphi$-functions
in~\cite{SchmelzerTrefethen2007}. The Pad{\'e} approximation works for general matrices $\bfL$, but the CF method used in \cite{SchmelzerTrefethen2007} works only if $z$ is negative and real, so that $\bfL$ must be symmetric negative definite. 

For large $N$ it will generally be too expensive to evaluate $\varphi_{\left[m,n\right]}(\bfL)$ in this form.
But suppose we have that $m \leq n$,
$q_n$ has $n$ distinct roots denoted by $s_1, \ldots, s_n$, 
and 
$p_m$ and $q_n$ have no roots in common. Then we can find a partial fraction expansion
$$
	\varphi_{\left[m,n\right]}(z) = c_0 + \sum_{i=1}^{n}\dfrac{c_k}{s_i - z},
$$
where $c_0$ is some constant and $c_k = \mathrm{Res}\left[\varphi_{\left[m,n\right]}(z),\rho_k\right]$. In practice one would find these coefficients simply by clearing the denominators.

The product of $\varphi(\bfL)$ with some vector $\bfw$ therefore is
\begin{equation}
	\varphi(\bfL)\,\bfw \approx \varphi_{\left[m,n\right]}(\bfL)\,\bfw = c_0\,\bfw + \sum_{i=1}^{n}c_k\left(s_i\bfI - \bfL\right)^{-1}\bfw.
	\label{eqn:phi_rational_approx_general}
\end{equation}

See \cite{SchmelzerTrefethen2007} for a discussion on how a common set of poles can be used for the evaluation of different $\varphi_j$. While this approach requires a higher degree $n$ in the rational approximation to achieve a given accuracy, in the use of exponential integration this would still lead to a more efficient method overall since the same computations can be used to evaluate different $\varphi_\ell$.


\subsection{Contour Integration}
\label{subsec:contour_integral}

The last approach we consider is based on the Cauchy integral formula
$$
	\phi(z) = \dfrac{1}{2\pi i}\int_{\Gamma}\,\dfrac{\phi(s)}{s-z}\,ds
$$
for a fixed value of $z$, where $\phi(z)$ is some arbitrary function and $\Gamma$ is a contour in the complex plane that encloses $z$ and is well-separated from $0$. This formula still holds when replacing $z$ by some general matrix $\bfL$, so that
\begin{equation}
	\phi(\bfL) = \dfrac{1}{2\pi i}\int_{\Gamma}\,\phi(s)\,\left(s\bfI-\bfL\right)^{-1}\,ds,
	\label{eqn:contour_integral_L}
\end{equation}
where $\Gamma$ can be any contour that encloses all the eigenvalues of $\bfL$. The integral is then approximated using some quadrature rule.

There is some freedom in choosing the contour integral and the quadrature rule. 
Kassam and Trefethen~\cite{KassamTrefethen2005} proposed 
the contour integral approach
to circumvent the cancellation error in~\eqref{eqn:phi_recurrence_relation}.
For convenience, they let $\Gamma$ simply be a circle in the complex plane that is large enough to enclose all the eigenvalues of $\bfL$, and then used the trapezoidal rule for the approximation. If $\bfL$ is real, then one can additionally simplify the calculations by considering only points on the upper half of a circle with its center on the real axis, and taking the real part of the result. Discretizing the contour using $M$ points $s_i$
and using the trapezoidal rule to evaluate \eqref{eqn:contour_integral_L}, we have
\begin{equation}
	\phi(\bfL) \approx \dfrac{1}{M}\,\Re\left(\sum\limits_{i=1}^M\,s_i\,\phi(s_i)\,\left(s_i\bfI-\bfL\right)^{-1}\right),
\label{eqn:phi_contour_intergral_discretized}
\end{equation}
where $M$ 
must be chosen large enough to give a good approximation. Then let $\phi = \varphi$ and multiply by $\bfw$ to get an approximation of the product $\varphi(\bfL)\bfw$. This is important in the context of multi-stage ETDRK.

Since the same contour integral will be used throughout the procedure, the quadrature points $s_i$ remain the same for all $\varphi$-functions and one can therefore use the same solutions $\bfv_i = \left(s_i\bfI-\bfL\right)^{-1}\bfw$ of the resolvent systems when computing the product of different $\varphi$-functions with some vector $\bfw$. 

A contour integral that specifically applies to $\varphi$-functions is
$$
	\varphi_\ell (z) = \dfrac{1}{2\pi i}\int_{\Gamma}\,\dfrac{e^s}{s^\ell }\dfrac{1}{s-z}\,ds.
$$
Again, different contour integrals and quadrature rules can be used. If $z$ is on the negative real line or close to it,
then we can use a Hankel contour (see also \cite{TrefethenWeidemanSchmelzer2006}). Letting $z = \bfL$ and using the trapezoidal rule, we get
\begin{equation}
	\varphi_\ell (\bfL) \approx \dfrac{1}{M}\,\Re\left(\sum\limits_{i=1}^M\,\dfrac{e^{s_i}}{s_i^\ell }\,\left(s_i\bfI-\bfL\right)^{-1}\right),
\label{eqn:varphi_contour_intergral_discretized}
\end{equation}
for quadrature points $s_i$ on $\Gamma$. This integral representation has the advantage that the integrand is exponentially decaying and therefore fewer quadrature points need to be used~\cite{SchmelzerTrefethen2007, HochbruckOstermann2010}. 

Incidentally, \eqref{eqn:phi_contour_intergral_discretized}, \eqref{eqn:varphi_contour_intergral_discretized} and \eqref{eqn:phi_rational_approx_general} all require an efficient procedure for solving linear systems of the form
$$
	(s_i\bfI - \bfL)\,\bfv_i = \bfw .
$$
This can be achieved, for instance, using
sparse direct solvers or preconditioned Krylov subspace methods. Solving $M$ different linear systems might seem prohibitive, but when using Krylov methods one has the advantage that $\mathcal{K}_M(\bfL,\bfw) = \mathcal{K}_M(s\bfI - \bfL,\bfw)$ for all $s \in \mathbb{C}^N$, so that the same Krylov subspace can theoretically be used for all $s_i$. The computation also allows for parallelization since each system can be solved independently. 



\section{The Discrete Adjoint Method}
\label{sec:discrete_adjoint_method}

To save on notation, we omit throughout the next five sections the subscript $\tau$ and write
the abstract representation~\eqref{eqn:ETDRK_T} of the discrete linear time-stepping system
as
\begin{equation*}
	\bft\left(\bfy,\bfm\right) = \zero .	
\end{equation*}
 It is now used in conjunction with the adjoint method to systematically find the procedures for computing the action of the sensitivity matrix 
$$
	\bfJ = \frac{\partialsp\bfd}{\partialsp\bfm} = \frac{\partialsp\bfd}{\partialsp\bfy}\hs \frac{\partialsp\bfy}{\partialsp\bfm}.
$$
We emphasize that here $\bfy = \bfy (\bfm )$ is a discrete vector which is the solution of the discrete forward problem.
To find an expression for $\frac{\partial\bfy}{\partial\bfm}$, we start by differentiating $\bft\left(\bfy(\bfm),\bfm\right) = \zero$ with respect to $\bfm$:
\begin{align*}
	\frac{\partial}{\partial \bfm}\bft\left(\bfy(\bfm),\bfm\right) = \frac{\partialsp\bft}{\partialsp\bfm} + 
	\frac{\partialsp\bft}{\partialsp\bfy}\frac{\partialsp\bfy}{\partialsp\bfm} = \zero_{N \times N_{\bfm}}.
\end{align*}
It follows that
\begin{equation}
\frac{\partialsp\bfy}{\partialsp\bfm} = 
-\left(\frac{\partialsp\bft}{\partialsp\bfy}\right)^{-1}\dfrac{\partialsp\bft}{\partialsp\bfm} ,
	\quad {\rm hence}\;\bfJ = -\dfrac{\partialsp\bfd}{\partialsp\bfy}\hs\left(\dfrac{\partialsp\bft}{\partialsp\bfy}\right)^{-1}\dfrac{\partialsp\bft}{\partialsp\bfm} .
  \label{eqn:J_and_Jt}
\end{equation}
It is neither desirable nor necessary to compute the Jacobian $\bfJ$ explicitly: 
what one really needs is to be able to quickly compute the action of $\bfJ$ and its transpose 
on appropriately-sized arbitrary vectors $\bfw$. Notice that the computation of the product of $\bfJ$ with some vector $\bfw$ of length $N_{\bfm}$ requires the solution of the \textit{linearized forward problem} $\frac{\partialsp\bft}{\partialsp\bfy}\bfv = \bfq$ (with $\bfq = \frac{\partialsp\bft}{\partialsp\bfm}\bfw$ here); we will discuss the structure of $\frac{\partialsp\bft}{\partialsp\bfy}$ in the following section. The computation of the product of $\bfJ^\top$ with some vector $\bfw$ of length $N$ requires the solution of the \textit{adjoint problem} 
$\left(\frac{\partialsp\bft}{\partialsp\bfy}\right)^{\top}\bflambda = \bftheta$, 
where $\bflambda$ is the \textit{adjoint solution} and $\bftheta$ is the \textit{adjoint source}. 
In this context $\bftheta = \left(\frac{\partialsp\bfd}{\partialsp\bfy}\right)^{\top}\bfw$.

Now, to compute the gradient of the misfit function, we use the chain rule to write
$$
	\nabla_\bfm\opM = \left(\dfrac{\partialsp\bfd}{\partialsp\bfm}\right)^\top\nabla_\bfd\opM = \bfJ^\top\nabla_\bfd\opM.
$$
The gradient $\nabla_\bfm\opM$ is therefore easily obtained by letting $\bfw = \nabla_\bfd\opM$ above. In addition to the forward solution $\bfy$,
one must compute the adjoint solution in order to get the gradient.
 
The gradient is used by all iterative gradient-based optimization methods,
including steepest descent, nonlinear conjugate gradient, quasi-Newton such as BFGS, Gauss-Newton and Levenberg-Marquardt. 
The latter two methods, in particular, require the computation of a solution to the linearized forward problem as well as an adjoint solution, in addition to the computation of $\bfy$.



\section{The derivatives of $\bft$}
\label{sec:derivatives_of_T}

We now focus on deriving expressions of the derivatives for $\dfrac{\partialsp\bft}{\partialsp\bfy}$ and $\dfrac{\partialsp\bft}{\partialsp\bfm}$ required by the sensitivity matrix. Recall equations \eqref{eqn:ETDRK_T} and \eqref{eqn:ETDRK_miscellaneous_definitions}, and let $\bfA_{k}^s$ be the $s \times s$ block matrix where the $(i,j)$th entry is $a_{ij}(\tau_k\bfL_k)$. To ease the notation further, let
\begin{subequations}
\begin{align}
t_{k,i} &= t_{k} + c_{i}\tau_{k}, \\
\underline{\bfy}_{k+1,i} &= e^{c_i\tau_k{\bfL}_{k}(\bfy_{k},\bfm)}\bfy_{k} + \tau_k\sum\limits^{i-1}_{j=1}a_{ij}(\tau_k{\bfL}_{k}(\bfy_{k},\bfm))\bfY_{k+1,j},
\end{align}
and\footnote{The Jacobian $\dfrac{\partialsp\bfn_{k,i}}{\partialsp\bfy}$ is taken with respect to the semi-discretized 
$\bfy(t)$.}
\begin{equation}
\dfrac{\partialsp\bfN_{k+1}}{\partialsp\bfy} = \mdiag{\dfrac{\partialsp\bfn_{k,1}}{\partialsp\bfy}, \cdots , \dfrac{\partialsp\bfn_{k,s}}{\partialsp\bfy}}.
\end{equation}
\end{subequations}
The solution procedure at the $(k+1)$th time-step is represented by
\begin{equation}
\bft_{k+1} = \begin{bmatrix} \bfY_{k+1}-\bfN_{k+1} \\
	-e^{\tau_k\bfL_{k}}\bfy_{k} -\bfB_{k+1,s}^\top\bfY_{k+1} + \bfy_{k+1}\end{bmatrix}.
\end{equation}
We will now take the derivative of $\bft$ with respect to $\bfy$ and $\bfm$ in turn.


\subsection{Computing $\dfrac{\partialsp\bft}{\partialsp\bfy}$}
\label{subsec:derivative_dTdy}

We compute the derivative of $\dfrac{\partialsp\bft}{\partialsp\bfy} = \bfT - \dfrac{\partialsp\bfn}{\partialsp\bfy}$ at the $(k+1)$th time-step. 
Letting $\widehat{\bfy}_{k+1} = \begin{bmatrix} \bfY_{k+1}^\top & \bfy_{k+1}^\top \end{bmatrix}^\top$, 
this derivative is
\begin{align*}
	\dfrac{\partialsp\bft_{k+1}}{\partialsp\bfy} &= \begin{bmatrix} \dfrac{\partialsp\bft_{k+1}}{\partialsp\widehat{\bfy}_1} & \dfrac{\partialsp\bft_{k+1}}{\partialsp\widehat{\bfy}_2} & \cdots & \dfrac{\partialsp\bft_{k+1}}{\partialsp\widehat{\bfy}_K}\end{bmatrix},
\end{align*}
with
\begin{equation}
\dfrac{\partialsp\bft_{k+1}}{\partialsp\widehat{\bfy}_j} = \begin{bmatrix} \dfrac{\partial}{\partialsp\widehat{\bfy}_j} \left(\bfY_{k+1}-\bfN_{k+1}\right) \\
	\dfrac{\partial}{\partialsp\widehat{\bfy}_j}\left(-e^{\tau_k\bfL_{k}}\bfy_{k} -\bfB_{k+1}^\top\bfY_{k+1} + \bfy_{k+1}\right)\end{bmatrix}.
	\label{eqn:ETDRK_dTkdYj_general}
\end{equation}
Looking at the terms in \eqref{eqn:ETDRK_dTkdYj_general} individually and using the chain rule, we have for $j=k+1$
$$
\begin{aligned}
	\dfrac{\partialsp\bfY_{k+1}}{\partialsp\widehat{\bfy}_{k+1}} &= \begin{bmatrix} \bfI_{sN\times sN} & \zero_{sN \times N} \end{bmatrix} & \dfrac{\partialsp\bfN_{k+1}}{\partialsp\widehat{\bfy}_{k+1}} &= \tau_k \dfrac{\partialsp\bfN_{k+1}}{\partialsp\bfy}\begin{bmatrix} \bfA_{k}^s & \zero_{sN\times N} \end{bmatrix} \\
	\dfrac{\partialsp\bfy_{k+1}}{\partialsp\widehat{\bfy}_{k+1}} &= \begin{bmatrix} \zero_{N \times sN} & \bfI_{N\times N} \end{bmatrix} & \dfrac{\partial(\bfB_{k+1}^\top\bfY_{k+1})}{\partialsp\widehat{\bfy}_{k+1}} &= \begin{bmatrix} \bfB_{k+1}^\top & \zero_{N \times N} \end{bmatrix}
\end{aligned}
$$
and for $j = k$
$$
\begin{aligned}
\dfrac{\partialsp\bfN_{k+1}}{\partialsp\widehat{\bfy}_{k}} &= \begin{bmatrix} \zero_{sN\times sN} & \dfrac{\partialsp\bfN_{k+1}}{\partialsp\bfy_k} \end{bmatrix}\quad & 
\dfrac{\partial(e^{\tau_k\bfL_{k}}\bfy_{k})}{\partialsp\widehat{\bfy}_{k}} &= \begin{bmatrix} \zero_{N \times sN} & \dfrac{\partial(e^{\tau_k\bfL_k}\bfy_k)}{\partialsp\bfy_k} \end{bmatrix} \\
 & & \dfrac{\partial(\bfB_{k+1}^\top\bfY_{k})}{\partialsp\widehat{\bfy}_k} &= \begin{bmatrix} \zero_{N \times sN} & \dfrac{\partial(\bfB_{k+1}^\top\bfY_{k+1})}{\partialsp\bfy_k} \end{bmatrix}.
\end{aligned}
$$
The terms
\begin{align*}
\dfrac{\partial(\bfB_{k+1}^\top\bfY_{k+1})}{\partialsp\bfy_k} &= \tau_k\sum_{i=1}^s\dfrac{\partial(b_{i}(\tau_k\bfL_k(\bfy_k))\bfY_{k+1,i})}{\partialsp\bfy_k} \\
\dfrac{\partialsp\bfn_{k,i}\left(\underline{\bfy}_{k,i}(\bfy_k),\bfy_k\right)}{\partialsp\bfy_k} &= \dfrac{\partialsp\bfn_{k,i}}{\partialsp\bfy}\dfrac{\partialsp\underline{\bfy}_{k+1,i}\left(\bfy_k\right)}{\partialsp\bfy_k} + \dfrac{\partialsp\bfn_{k,i}\left(\bfy_k\right)}{\partialsp\bfy_k}
\end{align*}
with
$$
\dfrac{\partialsp\underline{\bfy}_{k+1,i}\left(\bfy_k\right)}{\partialsp\bfy_k} = \dfrac{\partialsp e^{c_i\tau_k\bfL_k(\bfy_k)}\bfy_k}{\partialsp\bfy_k} + \tau_k\sum\limits^{i-1}_{j=1}\dfrac{\partial(a_{ij}(\tau_k\bfL_k(\bfy_k))\bfY_{k+1,j})}{\partialsp\bfy_k},
$$
where 
$$
\dfrac{\partialsp (e^{c_i\tau_k\bfL_k(\bfy_k)}\bfy_k)}{\partialsp\bfy_k} = e^{c_i\tau_k\bfL_k(\bfy_k)} + \dfrac{\partialsp (e^{c_i\tau_k\bfL_k(\bfy_k)}\bfy_k^{\mathrm{fixed}})}{\partialsp\bfy_k},
$$
will be discussed in detail in the following section.

Now let
\begin{equation}
\begin{aligned}
	\bfA_{k+1} &= \bfI_{sN} - \tau_k\hs \dfrac{\partialsp\bfN_{k+1}}{\partialsp\bfy}\hs\bfA_k^s \qquad \bfC_{k+1} = -\dfrac{\partialsp\bfN_{k+1}}{\partialsp\bfy_k}, \\
	\bfD_{k+1} &= -\dfrac{\partialsp\left(e^{\tau_k\bfL_k}\bfy_k\right)}{\partialsp\bfy_k} - \dfrac{\partialsp\left(\bfB_{k+1}^\top\bfY_{k+1}\right)}{\partialsp\bfy_k}
\label{eqn:ETDRK_independent_Ak_Ck_Dk_definitions}
\end{aligned}
\end{equation}
so that the $(k+1,j)$th $(s+1)N \times (s+1)N$ block of $\dfrac{\partialsp\bft}{\partialsp\bfy}$ is
\begin{equation}
	\dfrac{\partialsp\bft_{k+1}}{\partialsp\widehat{\bfy}_j} = \begin{cases}
	\begin{bmatrix} \bfA_{k+1} & \zero_{sN\times N} \\ -\bfB_{k+1}^\top & \bfI_N \end{bmatrix} & \quad\mathrm{if}\;j=k+1 \\
	\begin{bmatrix} \zero_{sN\times sN} & \bfC_{k+1} \\ \zero_{N\times sN}  & \bfD_{k+1} \end{bmatrix} & \quad\mathrm{if}\;j=k \\
	\quad\zero_{(s+1)N \times (s+1)N} & \quad\mathrm{otherwise}.
	\end{cases}
\end{equation}
and we can therefore write
\begin{equation}
\dfrac{\partialsp\bft}{\partialsp{\bfy}} = \begin{bmatrix}
   \bfA_1 \\
     -\bfB_{1}^\top  &  \bfI_N  \\
                & \bfC_2 & \bfA_2 \\
                & \bfD_2 & -\bfB_{2}^\top & \bfI_N \\
                &  &    \ddots      & \ddots & \ddots  \\
                &  & & &    \bfC_{K} & \bfA_K \\
           & & & &  \bfD_K & -\bfB_{K}^\top    & \bfI_N
\end{bmatrix},
\label{eqn:ETDRK_dTdY_general}
\end{equation}
which is a block lower-triangular matrix representing the linearized forward problem.


\subsection{Computing $\dfrac{\partialsp\bft}{\partialsp\bfm}$}

We now turn our attention to computing the derivative of $\bft$ with respect to $\bfm$. Consider the derivative of the $(k+1)$th time-step:
\begin{equation}
\dfrac{\partialsp\bft_{k+1}}{\partialsp\bfm} = \begin{bmatrix} -\dfrac{\partialsp\bfN_{k+1}(\bfy(\bfm),\bfm)}{\partialsp\bfm} \\
	-\dfrac{\partial}{\partialsp\bfm}\left(e^{\tau_k\bfL_k}\bfy_k - \bfB_{k+1}^\top\bfY_{k+1}\right) \end{bmatrix}.
	\label{eqn:ETD_dTkdM_setup}
\end{equation}
The individual terms are
\begin{subequations}
\begin{equation}
\dfrac{\partial(\bfB_{k+1}^\top\bfY_{k})}{\partialsp\bfm} = \tau_k\sum_{i=1}^s\dfrac{\partial(b_{i}(\tau_k\bfL_k(\bfm))\bfY_{k+1,i})}{\partialsp\bfm},
\end{equation}
and
\begin{equation}
\dfrac{\partialsp\bfn_{k,i}\left(\underline{\bfy}_{k+1,i}(\bfm),\bfm\right)}{\partialsp\bfm} = \dfrac{\partialsp\bfn_{k,i}(\bfm)}{\partialsp\bfy}\dfrac{\partialsp\underline{\bfy}_{k+1,i}}{\partialsp\bfm} + \dfrac{\partialsp\bfn_{k,i}(\bfm)}{\partialsp\bfm}
\end{equation}
with
$$
	\dfrac{\partialsp\underline{\bfy}_{k+1,i}}{\partialsp\bfm} = \dfrac{\partial\left(e^{c_i\tau_k\bfL_k(\bfm)}\bfy_k\right)}{\partialsp\bfm} + \tau_k\sum\limits^{i-1}_{j=1}\dfrac{\partial(a_{ij}(\tau_k\bfL_k(\bfm))\bfY_{k+1,i})}{\partialsp\bfm}.
$$
\end{subequations}
If the linear terms $\tau_k\bfL_k$ are independent of $\bfm$ then the only term that depends on $\bfm$ is $\bfN_k$, so trivially we have
\begin{equation}
\dfrac{\partialsp\bft_{k+1}}{\partialsp\bfm} = \begin{bmatrix} -\dfrac{\partialsp\bfN_{k+1}(\bfm)}{\partialsp\bfm} \\
	\zero_{N \times N_{\bfm}}\end{bmatrix}.
	\label{eqn:ETD_independent_L_dTkdM_setup}
\end{equation}
We consider the case of $\tau_k\bfL_k$ depending on $\bfm$ in the next section.

The product of $\dfrac{\partialsp\bft_{k+1}}{\partialsp\bfm}$ with some vector $\bfw_{\bfm}$ of length $N_{\bfm}$ is obvious. The product $\dfrac{\partialsp\bft_{k+1}}{\partialsp\bfm}^\top\widehat{\bfw}_{k+1}$, with $\widehat{\bfw}_{k+1} = \begin{bmatrix} \bfW_{k+1}^\top & \bfw_{k+1}^\top \end{bmatrix}^\top$ an arbitrary vector of length $(s+1)N$ defined analogously to $\widehat{\bfy}_{k+1}$ (replacing $\bfY_{k+1}$ by $\bfW_{k+1}$), is
\begin{eqnarray}
\dfrac{\partialsp\bft_{k+1}}{\partialsp\bfm}^\top\widehat{\bfw}_{k+1} &=& 
-\sum\limits_{i=1}^s\dfrac{\partialsp\bfn_{k,i}}{\partialsp\bfm}^\top\bfW_{k+1,i} - \left\{\dfrac{\partialsp(e^{\tau_k\bfL_k(\bfm)}\bfy_k)}{\partialsp\bfm}^\top\bfw_{k+1} + \right. \nonumber \\
&+& \left. \sum\limits_{i=1}^s\left(\dfrac{\partialsp\underline{\bfy}_{k+1,i}}{\partialsp\bfm}^\top\bfv_{k+1,i} + \tau_k{\dfrac{\partial(b_{i}(\tau_k\bfL_k(\bfm))\bfY_{k+1,i})}{\partialsp\bfm}}^\top\bfw_{k+1}\right)\right\},
\label{eqn:ETD_dTkdM_transpose_times_w}
\end{eqnarray}
where $\bfv_{k+1,i} = \dfrac{\partialsp\bfn_{k,i}}{\partialsp\bfy}^\top\bfW_{k+1,i}$ and
$$
\dfrac{\partialsp\underline{\bfy}_{k+1,i}}{\partialsp\bfm}^\top\bfv_{k+1,i} = {\dfrac{\partialsp (e^{c_i\tau_k\bfL_k(\bfm)}\bfy_k)}{\partialsp\bfm}}^\top \bfv_{k+1,i} + \tau_{k}\sum\limits^{i-1}_{j=1}{\dfrac{\partial(a_{ij}(\tau_k\bfL_k(\bfm))\bfY_{k+1,i})}{\partialsp\bfm}}^\top\bfv_{k+1,i}.
$$
The terms in braces in \eqref{eqn:ETD_dTkdM_transpose_times_w} are ignored if $\bfL_k$ is independent of $\bfm$.



\section{The Derivatives of $\varphi$}
\label{sec:derivatives_of_phi}

If the linear operator $\bfL_{k}$ depends on either $\bfy_{k}$ or $\bfm$ we have to be able to take the derivatives of the products $a_{ij}(\tau_k\bfL_{k}(\bfy_{k},\bfm))\bfw$, $b_{i}(\tau_k\bfL_{k}(\bfy_{k},\bfm))\bfw$ and $e^{c_i\tau_k\bfL_{k}(\bfy_{k},\bfm)}\bfw$, where $\bfw$ is an arbitrary vector of length $N$, with respect to $\bfy_{k}$ or $\bfm$.

Since the $a_{ij}$ and $b_i$ are linear combinations of the $\varphi_{\ell,i}$-and $\varphi_\ell$-functions respectively, we will simply consider the derivatives of $\varphi(\bfL)\bfw$ for some arbitrary $\varphi$-function and matrix $\bfL = \tau_k{\bfL}_{k}(\bfy_{k},\bfm)$. Without loss of generality we assume in this section that 
%
the differentiation is with respect to $\bfy_{k}$, but all the results also apply when taking the derivative with respect to $\bfm$.

The calculation of the derivatives of the products of the $\varphi_\ell$-functions will depend on the way these terms are evaluated numerically. Recall the approaches for evaluating $\varphi_\ell$ reviewed in Section~\ref{sec:action_of_phi}.


\subsection{Krylov Subspace Methods}

Using Krylov subspace methods is unfortunately unsuitable for our purposes since if $\bfL$ depends on $\bfy_k$ or $\bfm$, it is extremely difficult to find the dependence of the right-hand side of \eqref{eqn:phi_product_Krylov} on these variables.
For this reason we will not consider this approach further here, although it can of course be used for parameter estimation problems where ${\bfL}$ is independent of $\bfy_k$ and $\bfm$.


\subsection{Polynomial Approach}

If we represent $\varphi$ by a truncated polynomial series as in \eqref{eqn:phi_product_general_polynomial} and multiply by $\bfw$, we have
$$
	\varphi(\bfL(\bfy_k))\bfw \approx \sum_{j=0}^M\,c_j\,\bfL^j(\bfy_k)\bfw.
$$
Let $\bfL$ depend on a single parameter $z$ first. The derivative of $\bfL^j(z)\bfw$ then is
$$
	\dfrac{\partial\,\bfL^j(z)}{\partial\,z}\bfw = \sum\limits_{i=1}^j\bfL^{i-1}\dfrac{\partial\,\bfL(z)}{\partial\,z}\bfL^{j-i}\bfw = \sum\limits_{i=1}^j\bfL^{i-1}\dfrac{\partial\,\bfL(z)}{\partial\,z}\bfv_{j-i}
$$
with $\bfv_{j-i} = \bfL^{j-i}\bfw$. In the case of $\bfL = \bfL(\bfy_k)$ we therefore have
$$
	\dfrac{\partial\,\bfL^j(\bfy_k)\bfw}{\partial\,\bfy_k} = \sum\limits_{i=1}^j\bfL^{i-1}\dfrac{\partial\,\bfL(\bfy_k)\bfv_{j-i}}{\partial\,\bfy_k}.
$$
This implies that we need to have $j-1$ derivatives $\dfrac{\partial\,\bfL(\bfy_k)\bfv_{j-i}}{\partial\,\bfy_k}$ available for each $j$, but of course a lot of these derivatives can also be reused for different values of $j$.

In the special case where the matrices $\bfL$ and $\dfrac{\partial\,\bfL}{\partial\,y}$ commute for each element $y$ of $\bfy_k$, we can simplify the above by using the matrix analogue of the usual power rule:
$$
	\dfrac{\partial\,\bfL^j(\bfy_k)}{\partial\,\bfy_k}\bfw = \bfL^{j-1}\dfrac{\partial\,\bfL(\bfy_k)\bfw}{\partial\,\bfy_k}.
$$
This is much less cumbersome to implement, but it is hard to think of realistic scenarios where this commutativity property would hold.

The drawback of considering each monomial term on its own is that this does not necessarily reflect how the polynomial approximation of $\varphi$ is actually computed. For instance, the Chebyshev approximation \eqref{eqn:phi_product_chebyshev} is
$$
	\varphi(\bfL)\bfw \approx \sum_{j=0}^M\,c_j\,T_j(\bfL)\bfw
$$
assuming $\bfL(\bfm)$ is Hermitian or skew-Hermitian and the eigenvalues of $\bfL(\bfy_k)$ all lie inside $[-1,1]$. Each $T_j(\bfy_k)\bfw$ is computed using the recurrence relation
$$
	T_{j+1}(\bfL(\bfy_k))\bfw = 2\bfL\,T_{j}(\bfL(\bfy_k))\bfw - T_{j-1}(\bfL(\bfy_k))\bfw, \qquad j = 1,2,\ldots
$$
initialized by $T_{0}(\bfL(\bfy_k))\bfw = \bfw$ and $T_{1}(\bfL(\bfy_k))\bfw = \bfL(\bfy_k)\bfw$. Taking the derivative of \eqref{eqn:phi_product_chebyshev} with respect to $\bfy_k$,
$$
	\dfrac{\partial\,\varphi(\bfL)\bfw}{\partial\,\bfy_k} \approx \sum_{j=0}^M\,c_j\dfrac{\partial\,T_j(\bfL)\bfw}{\partial\,\bfy_k},
$$
with the recurrence relation
$$
	\dfrac{\partial\,T_{j+1}(\bfL)\bfw}{\partial\,\bfy_k} = 2\dfrac{\partial\,(\bfL\,\bfv)}{\partial\,\bfy_k} + 2\bfL\dfrac{\partial\,T_{j}(\bfL)\bfw}{\partial\,\bfy_k} - \dfrac{\partial\,T_{j-1}(\bfL)\bfw}{\partial\,\bfy_k}, \qquad j = 1,2,\ldots
$$
where $\bfv = T_{j}\bfw$, $\dfrac{\partial\,T_{0}(\bfL)\bfw}{\partial\,\bfy_k} = \zero_{N \times N}$ and $\dfrac{\partial\,T_{1}(\bfL)\bfw}{\partial\,\bfy_k} = \dfrac{\partial\,(\bfL\,\bfw)}{\partial\,\bfy_k}$.

We also need to be able to compute the products of the transposes of $\dfrac{\partial\,\varphi(\bfL)\bfw}{\partial\,\bfy_k}$ with some vector; 
this is straightforward to derive from the equations above.


\subsection{Rational Approximations and Contour Integration}

We saw in Sections \ref{subsec:rational_approximation} and \ref{subsec:contour_integral} that the action of a $\varphi$-function can be computed by both rational approximations (under certain conditions) and contour integrals in the form
$$
	\varphi(\bfL)\bfw \approx \sum\limits_{i=1}^M\,c_i\left(s_i\bfI-\bfL\right)^{-1}\bfw,
$$
for some complex scalars $c_i$ and $s_i$, where the $s_i$ do not coincide with the eigenvalues of $\bfL$.

To find the derivative of $\varphi(\bfL)\bfw$ in this case, let $\bfv_i = \left(s_i\bfI-\bfL\right)^{-1}\bfw$, so that
\begin{equation}
	\dfrac{\partial\,\varphi(\bfL)\,\bfw}{\partial\,\bfy_k} \approx \sum\limits_{i=1}^M\,c_i\dfrac{\partial\,\bfv_i}{\partial\,\bfy_k}.	\label{eqn:dphidM_times_w_general_partial_fraction}
\end{equation}
To find an expression for $\dfrac{\partial\,\bfv_i}{\partial\,\bfy_k}$, consider $\left(s_i\bfI-\bfL\right)\bfv_i = \bfw$ and take the derivative with respect to $\bfy_k$ on both sides,
\begin{align*}
\dfrac{\partial}{\partial\,\bfy_k}\left(s_i\bfI-\bfL\right)\bfv_i = \zero_{N\times N} & \qquad \Rightarrow \qquad \left(s_i\bfI-\bfL\right)\dfrac{\partial\,\bfv_i}{\partial\,\bfy_k} - \dfrac{\partial\,(\bfL\,\bfv_i)}{\partial\,\bfy_k} = \zero_{N\times N} \\
& \qquad \Rightarrow \qquad \dfrac{\partial\,\bfv_i}{\partial\,\bfy_k} = \left(s_i\bfI-\bfL\right)^{-1}\dfrac{\partial\,(\bfL\,\bfv_i)}{\partial\,\bfy_k},
\end{align*}
where the $\bfv_i$ on the right-hand side is taken to be fixed.

For each term in the sum in \eqref{eqn:dphidM_times_w_general_partial_fraction} we thus require two matrix solves, one to find $\bfv_i$ and an additional one to then find $\left(s_i\bfI-\bfL\right)^{-1}\frac{\partial\,(\bfL\bfv_i)}{\partial\,\bfy_k}$. As mentioned previously, the $\bfv_i$ can be computed using the same Krylov subspace. However, with $\bfz$ an arbitrary vector of length $N$, the linear systems
$$
	\left(s_i\bfI-\bfL\right)\bfu_i = \dfrac{\partial\,(\bfL\,\bfv_i)}{\partial\,\bfy_k}\bfz
$$
each have a different right-hand side, and therefore we are no longer working in the same Krylov subspace. This means that we need to solve $M$ different linear systems just for a single evaluation of the derivative of a given $\varphi$-function, which is not ideal. The process is fortunately highly parallelizable and given the ease of access to a large number of processors these days, we do not consider this to be the bottleneck it might have been just a few years ago. Nonetheless it is a significant inconvenience 
for many a mathematician, 
and the polynomial approach does not suffer from this limitation.

We also should be able to compute the transpose of \eqref{eqn:dphidM_times_w_general_partial_fraction}. In this case we have
$$
	\dfrac{\partial\,\varphi(\bfL)\,\bfw}{\partial\,\bfy_k}^\top \approx \sum\limits_{i=1}^M\,c_i\dfrac{\partial\,\bfv_i}{\partial\,\bfy_k}^\top,
$$
with $\dfrac{\partial\,\bfv_i}{\partial\,\bfy_k}^\top = \dfrac{\partial\,(\bfL\,\bfv_i)}{\partial\,\bfy_k}^\top\left(s_i\bfI-\bfL\right)^{-\top} = \dfrac{\partial\,(\bfL\,\bfv_i)}{\partial\,\bfy_k}^\top\left(s_i\bfI-\bfL^\top\right)^{-1}$, so when multiplying by some vector $\bfz$ we can work within the same Krylov subspace $\mathcal{K}_M(\bfL^\top,\bfz)$ for any 
$s_i \in \mathbb{C}^N$.



\section{Solving the Linearized Forward Problem}
\label{sec:linearized_forward_problem}

In Section \ref{sec:discrete_adjoint_method} we showed that computing the action of the sensitivity matrix $\bfJ$ on some vector $\bfw$ of length $N_{\bfm}$ requires the solution of the linearized forward problem
\begin{equation}
	\bfv = \left( \dfrac{\partialsp\bft}{\partialsp\bfy}\right)^{-1}\hs\bfq,	\label{eqn:linearized_forward_system_general}
\end{equation}
where $\bfv = \begin{bmatrix} \bfV_1^\top & \bfv_1^\top & \cdots & \bfV_K^\top & \bfv_K^\top \end{bmatrix}$ is the solution (the $\bfV_{k+1}$ represent the internal stages) and $\bfq = \begin{bmatrix} \bfQ_1^\top & \bfq_1^\top & \cdots & \bfQ_K^\top & \bfq_K^\top \end{bmatrix}$ is taken in this section to be some arbitrary source term including internal stages. In the context of sensitivity analysis we will have $\bfq = \frac{\partialsp\bft}{\partialsp\bfm}\bfw$, 
see~\eqref{eqn:J_and_Jt}, where $\bfw$ is an arbitrary vector of length $N_{\bfm}$. The linearized time-stepping system $\frac{\partialsp\bft}{\partialsp\bfy}$ was derived in 
Section \ref{subsec:derivative_dTdy} and is given in \eqref{eqn:ETDRK_dTdY_general}. 

\subsection{The General Linearized ETDRK Method}
\label{subsec:general_linearized_forward_procedure}

The linear system to be solved is 
\begin{equation}
\begin{bmatrix}
   \bfA_1 \\
     -\bfB_{1}^\top  &  \bfI_N  \\
                & \bfC_2 & \bfA_2 \\
                & \bfD_2 & -\bfB_{2}^\top & \bfI_N \\
                &  &    \ddots      & \ddots & \ddots  \\
                &  & & &    \bfC_K & \bfA_K \\
           & & & &  \bfD_K & -\bfB_K^\top    & \bfI_N
\end{bmatrix}\begin{bmatrix}
   \bfV_1 \\
   \bfv_1 \\
   \bfV_2 \\
   \vdots \\
   \bfv_{K-1} \\
   \bfV_K \\
   \bfv_K
\end{bmatrix} = \begin{bmatrix}
   \bfQ_1 \\
   \bfq_1 \\
   \bfQ_2 \\
   \vdots \\
   \bfq_{K-1} \\
   \bfQ_K \\
   \bfq_K
\end{bmatrix},	\label{eqn:ETDRK_linearized_forward_system}
\end{equation}
with $\bfA_{k+1}$, $\bfC_{k+1}$ and $\bfD_{k+1}$ defined in \eqref{eqn:ETDRK_independent_Ak_Ck_Dk_definitions}
and $\bfB_{k+1}^{\top} = \tau_{k}\begin{bmatrix} b_1(\tau_{k}\bfL_{k}) & \cdots & b_s(\tau_{k}\bfL_{k}) \end{bmatrix}$. 
Since the system is block-lower triangular we use forward substitution to get
$$
	\bfv_{k+1} = \bfq_{k+1} + \bfB_{k+1}^{\top}\bfV_{k+1} - \bfD_{k+1}\bfv_{k}
$$
for $k = 0,\ldots,K-1$ and internal stages
$$
\bfA_{k+1}\bfV_{k+1} = \bfQ_{k+1} - \bfC_{k+1}\bfv_{k},
$$
with $\bfv_0 = \zero_{N \times 1}$.

Substituting \eqref{eqn:ETDRK_independent_Ak_Ck_Dk_definitions}
then gives
$$
	\bfv_{k+1} = \bfq_{k+1} + \bfB_{k+1}^{\top}\bfV_{k+1} + \dfrac{\partialsp\left(e^{\tau_{k}\bfL_{k}}\bfy_{k}\right)}{\partialsp\bfy_{k}}\bfv_{k} + \dfrac{\partialsp\left(\bfB_{k+1}^\top\bfY_{k+1}\right)}{\partialsp\bfy_{k}}\bfv_{k}
$$
for $k = 0,\ldots,K-1$ and internal stages
$$
\bfV_{k+1} = \bfQ_{k+1} + \dfrac{\partialsp\bfN_{k+1}}{\partialsp\bfy_{k}}\bfv_{k} + \tau_k\hs\dfrac{\partialsp\bfN_{k+1}}{\partialsp\bfy}\hs \bfA_{k}^s\bfV_{k+1}.
$$
The detailed solution procedure is summarized in Algorithm \ref{alg:ETDRK_linearized_forward_solution}.

\begin{algorithm}
\label{alg:ETDRK_linearized_forward_solution}
{\rm The Linearized ETDRK Method} \\[0.5em]
Let $\bfv_{0} = \zero_{N \times 1}$ and ignore the terms in braces if $\bfL_{k}$ does not depend on $\bfy_{k}$. \\
For $k = 0,\ldots,K-1$:
\begin{itemize}
\item For $i = 1,\ldots,s$, compute the internal stages 
\begin{subequations}
\begin{equation}
\begin{aligned}
\bfV_{k+1,i} &= \bfQ_{k+1,i} + \dfrac{\partialsp\bfn_{k,i}}{\partialsp\bfy}\left(e^{c_i\tau_k\bfL_k}\bfv_k + \tau_k\sum\limits_{j=1}^{i-1}\hs a_{ij}(\tau_k\bfL_k)\bfV_{k+1,j}\right) + \\
& \quad + \left\{\dfrac{\partialsp\bfn_{k,i}}{\partialsp\bfy_k}\bfv_k + \dfrac{\partialsp\bfn_{k,i}}{\partialsp\bfy}\left(\dfrac{\partialsp(e^{c_i\tau_k\bfL_k(\bfy_k)}\bfy_k^{\mathrm{fixed}})}{\partialsp\bfy_k}\bfv_k\right.\right. + \\
& \qquad\qquad + \left.\left.\tau_k\sum\limits^{i-1}_{j=1}\dfrac{\partial(a_{ij}(\tau_k\bfL_k(\bfy_k))\bfY_{k+1,j})}{\partialsp\bfy_k}\bfv_k\right)\right\}.
\end{aligned}
\label{eqn:ETDRK_linearized_forward_internal_stages}
\end{equation}
\item Compute the update:
\begin{equation}
\begin{aligned}
	\bfv_{k+1} &= \bfq_{k+1} + e^{\tau_k\bfL_k}\bfv_k + \tau_k\sum\limits^{s}_{i=1}\hs b_{i}(\tau_k\bfL_k)\bfV_{k+1,i} + \\
	& \quad + \left\{\dfrac{\partialsp (e^{\tau_k\bfL_k(\bfy_k)}\bfy_k^{\mathrm{fixed}})}{\partialsp\bfy_k}\bfv_k + \tau_k\sum_{i=1}^s\dfrac{\partial(b_{i}(\tau_k\bfL_k(\bfy_k))\bfY_{k+1,i})}{\partialsp\bfy_k}\bfv_k\right\}.
\end{aligned}
\label{eqn:ETDRK_linearized_forward_update}
\end{equation}
\label{eqn:ETDRK_linearized_forward_algorithm}
\end{subequations}
\end{itemize} \qed
\end{algorithm}

The derivatives in the braces in \eqref{eqn:ETDRK_linearized_forward_algorithm} can be evaluated by writing $a_{ij}$ and $b_{i}$ in terms of $\varphi$-functions and then differentiating the $\varphi$-functions as in Section \ref{sec:derivatives_of_phi}. Notice that without the terms in braces we simply have the standard ETDRK method where the nonlinear term has been linearized, except that we allow for the possibility of a source term in the update formula.


\subsection{Application to Krogstad's scheme}
\label{sec7.2}

We show 
the above algorithm 
applied to the scheme proposed in~\cite{Krogstad2005}.
Recall that the matrix coefficients $a_{ij}$ and $b_{i}$ are given in~\eqref{eqn:ETDRK_Butcher_tableau_Krogstad} in terms of the $\varphi$-functions. The linearized scheme is presented in Algorithm \ref{alg:ETDRK_linearized_forward_{k+1}rogstad_scheme}.

\begin{algorithm}
\label{alg:ETDRK_linearized_forward_{k+1}rogstad_scheme}
{\rm The Linearized Scheme} \\[0.5em]
Let $\bfv_{0} = \zero_{N \times 1}$, $\varphi_{\ell} = \varphi_{\ell}(\tau_k\bfL_k(\bfy_k))$ and 
$\varphi_{\ell,i} = \varphi_{\ell,i}(\tau_k\bfL_k(\bfy_k))$, and ignore the terms in braces if $\bfL_{k}$ does not depend on $\bfy_{k}$. For $k = 0,\ldots,K-1$:
\begin{itemize}
\item 
Compute the internal stages 
\begin{align*}
\bfV_{k+1,1} &= \bfQ_{k+1,1} + \dfrac{\partialsp\bfn_{k,1}}{\partialsp\bfy}\bfv_k + \left\{\dfrac{\partialsp\bfn_{k,1}}{\partialsp\bfy_k}\bfv_k\right\} \\
\bfV_{k+1,2} &= \bfQ_{k+1,2} + \dfrac{\partialsp\bfn_{k,2}}{\partialsp\bfy}\left(\varphi_{0,2}\bfv_k + \frac{\tau_k}{2}\varphi_{1,2}\bfV_{k+1,1}\right) + \\
& \quad + \left\{\dfrac{\partialsp\bfn_{k,2}}{\partialsp\bfy}\left(\dfrac{\partialsp (\varphi_{0,2}\hs \bfy_k^{\mathrm{fixed}})}{\partialsp\bfy_k} + \frac{\tau_k}{2} \dfrac{\partial(\varphi_{1,2}\hs \bfY_{k+1,1})}{\partialsp\bfy_k}\right)\bfv_k + \dfrac{\partialsp\bfn_{k,2}}{\partialsp\bfy_k}\bfv_k\right\}
\end{align*}
\begin{align*}
\bfV_{k+1,3} &= \bfQ_{k+1,3} + \dfrac{\partialsp\bfn_{k,3}}{\partialsp\bfy}\left(\varphi_{0,3}\bfv_k + \frac{\tau_k}{2}\varphi_{1,3}\bfV_{k+1,1} + \tau_k\varphi_{2,3}\left(\bfV_{k+1,2} - \bfV_{k+1,1}\right)\right) + \\
& \quad + \left\{\dfrac{\partialsp\bfn_{k,3}}{\partialsp\bfy_k}\bfv_k + \dfrac{\partialsp\bfn_{k,3}}{\partialsp\bfy}\left(\dfrac{\partial(\varphi_{0,3}\bfy_k^{\mathrm{fixed}})}{\partialsp\bfy_k} + \dfrac{\tau_k}{2}\dfrac{\partial(\varphi_{1,3}\bfY_{k+1,1})}{\partialsp\bfy_k}\right)\bfv_k \right\} \\
& \quad + \left\{\tau_k\dfrac{\partialsp\bfn_{k,3}}{\partialsp\bfy}\left(\dfrac{\partial(\varphi_{2,3}(\bfY_{k+1,2} - \bfY_{k+1,1}))}{\partialsp\bfy_k}\right)\bfv_k\right\} \\
\bfV_{k+1,4} &= \bfQ_{k+1,4} + \dfrac{\partialsp\bfn_{k,4}}{\partialsp\bfy}\left(\varphi_{0}\bfv_k + \tau_k\varphi_{1,4}\bfV_{k+1,1} + 2\tau_k\varphi_{2,4}\left(\bfV_{k+1,3}-\bfV_{k+1,1}\right)\right) + \\
& \quad + \left\{\dfrac{\partialsp\bfn_{k,4}}{\partialsp\bfy_k}\bfv_k + \dfrac{\partialsp\bfn_{k,4}}{\partialsp\bfy}\left(\dfrac{\partialsp(\varphi_{0}\bfy_k^{\mathrm{fixed}})}{\partialsp\bfy_k} + \tau_k\dfrac{\partial(\varphi_{1,4}\bfY_{k+1,1})}{\partialsp\bfy_k}\right)\bfv_k \right\} + \\
& \quad + \left\{2\tau_k\dfrac{\partialsp\bfn_{k,4}}{\partialsp\bfy}\left(\dfrac{\partial(\varphi_{2,4}(\bfY_{k+1,3} - \bfY_{k+1,1}))}{\partialsp\bfy_k}\right)\bfv_k\right\}
\end{align*}
\item Compute the update:
\begin{align*}
	\bfv_{k+1} &= \bfq_{k+1} + \varphi_{0}\bfv_k + \tau_k\varphi_{2}(-3\bfV_{k+1,1}+2\bfV_{k+1,2}+2\bfV_{k+1,3}-\bfV_{k+1,4}) + \\
	 & \qquad\qquad + \tau_k\varphi_{1}\bfV_{k+1,1} + 4\tau_k\varphi_{3}(\bfV_{k+1,1}-\bfV_{k+1,2}-\bfV_{k+1,3}+\bfV_{k+1,4}) + \\
	& \qquad\qquad + \left\{\dfrac{\partialsp (\varphi_{0}\bfy_k^{\mathrm{fixed}})}{\partialsp\bfy_k}\bfv_k + \tau_k\dfrac{\partial(\varphi_{1}\bfY_{k+1,1})}{\partialsp\bfy_k}\bfv_k + \right. \\
	& \qquad\qquad\qquad + \tau_k\dfrac{\partial(\varphi_{2}(-3\bfY_{k+1,1}+2\bfY_{k+1,2}+2\bfY_{k+1,3}-\bfY_{k+1,4}))}{\partialsp\bfy_k}\bfv_k +  \\
	& \qquad\qquad\qquad\qquad + \left. 4\tau_k\dfrac{\partial(\varphi_{3}(\bfY_{k+1,1}-\bfY_{k+1,2}-\bfY_{k+1,3}+\bfY_{k+1,4}))}{\partialsp\bfy_k}\bfv_k\right\}.
\end{align*}
\end{itemize} \qed \\[0.7em]
\end{algorithm}
The way the scheme is written here may of course not represent the most computationally efficient implementation: 
for instance, there are some opportunities for parallelization and one should also exploit the fact that $c_2 = c_3 = 1/2$.



\section{Solving the Adjoint Problem}
\label{sec:adjoint_problem}

The major ingredient needed when calculating the action of the transpose of the sensitivity matrix, and therefore especially important for the gradient computation, is the solution of the adjoint problem,
\begin{equation}
	\bflambda = \left( \dfrac{\partialsp\bft}{\partialsp\bfy}\right)^{-\top}\hs\bftheta,
	\label{eqn:adjoint_system_general}
\end{equation}
where $\bflambda = \begin{bmatrix} \bfLambda_1^\top & \bflambda_1^\top & \cdots & \bfLambda_K^\top & \bflambda_K^\top \end{bmatrix}^\top$ is the adjoint solution (with internal stages $\bfLambda_k$) and $\bftheta = \begin{bmatrix} \bfTheta_1^\top & \bftheta_1^\top & \cdots & \bfTheta_K^\top & \bftheta_K^\top \end{bmatrix}$ is taken in this section to be some arbitrary adjoint source; for the gradient computation we will have $\bftheta = \frac{\partialsp\bfd}{\partialsp\bfy}^\top\hs\nabla_\bfd\mathcal{M}$. 


\subsection{The General Adjoint ETDRK Method}
\label{subsec:general_adjoint_procedure}

The linearized time-stepping system $\frac{\partialsp\bft}{\partialsp\bfy}$ is lower block-triangular (see \eqref{eqn:ETDRK_dTdY_general}) and therefore the linearized forward problem can conceptually be seen as being solved by forward substitution, which corresponds to solving the problem forward in time, as we saw in the previous section. The adjoint, i.e., its transpose in the finite-dimensional setting we are in, will thus be upper block-triangular:
\begin{equation}
	\begin{bmatrix}
   \bfA_1^\top & -\bfB_{1} \\
               &  \bfI_N & \bfC_2^\top & \bfD_2^\top \\
               &  & \bfA_2^\top & -\bfB_{2} \\
               &  &  & \ddots   & \ddots & \ddots  \\
                &  &  &  & \bfI_N & \bfC_{K}^\top & \bfD_{K}^\top \\
                &  &  &  &  & \bfA_K^\top & -\bfB_{K} \\
           & & & &  &  & \bfI_N
\end{bmatrix}\hs \begin{bmatrix}
   \bfLambda_1 \\
   \bflambda_1 \\
   \bfLambda_2 \\
   \vdots \\
   \bflambda_{K-1} \\
   \bfLambda_K \\
   \bflambda_K
\end{bmatrix} = \begin{bmatrix}
   \bfTheta_1 \\
   \bftheta_1 \\
   \bfTheta_2 \\
   \vdots \\
   \bftheta_{K-1} \\
   \bfTheta_K \\
   \bftheta_K
\end{bmatrix},
	\label{eqn:ETD_adjoint_system}
\end{equation}
which is solved by backward substitution (this corresponds to solving the adjoint problem backward in time), 
for $k = K,K-1,\ldots,1$. One such step reads
\begin{subequations}
\begin{equation}
\bflambda_{k} = \bftheta_{k} - \bfC_{k+1}^\top\bfLambda_{k+1} - \bfD_{k+1}^\top\bflambda_{k+1}
\end{equation}
followed by the internal stages
\begin{equation}
\bfA_k^\top\bfLambda_{k} = \bfTheta_{k}  +\bfB_{k}\bflambda_{k},
\end{equation}
\end{subequations}
where $\bflambda_{K+1} = \zero_{N \times 1}$ and $\bfLambda_{K+1} = \zero_{sN \times 1}$.

Using the general formulas for $\bfA_k$, $\bfC_k$ and $\bfD_k$ defined in  \eqref{eqn:ETDRK_independent_Ak_Ck_Dk_definitions}
gives
\begin{subequations}
\begin{equation}
\bflambda_{k} = \bftheta_{k} + \left(\dfrac{\partialsp\bfN_{k+1}}{\partialsp\bfy}\right)^\top\bfLambda_{k+1} + \left(\dfrac{\partialsp (e^{\tau_k{\bfL}_{k}}\bfy_{k})}{\partialsp\bfy_{k}}^\top + \dfrac{\partialsp (\bfB_{k+1}^\top\bfY_{k+1})}{\partialsp\bfy_{k}}^\top\right)\bflambda_{k+1}
\end{equation}
with internal stages
\begin{equation}
\bfLambda_{k} = \bfTheta_{k} + \tau_{k-1}\left(\bfA_{k-1}^s\right)^\top\left(\dfrac{\partialsp\bfN_{k}}{\partialsp\bfy}\right)^\top\bfLambda_{k} + \bfB_{k}\bflambda_{k}.
\end{equation}
\end{subequations}
The detailed solution procedure is summarized in Algorithm \ref{alg:ETDRK_adjoint_solution}.

\begin{algorithm}
\label{alg:ETDRK_adjoint_solution}
{\rm The Adjoint Exponential Runge-Kutta Time-Stepping Method} \\[0.5em]
For $k = K, K-1,\ldots,1$:
\begin{itemize}
\item Compute
\begin{subequations}
\begin{eqnarray}
\bflambda_{k} &=& \bftheta_{k} + e^{\tau_k\bfL_{k}^\top}\bflambda_{k+1} + 
\sum\limits_{i=1}^s e^{c_{i}\hs\tau_k\bfL_{k}^\top}\left(\dfrac{\partialsp\bfn_{k,i}}{\partialsp\bfy}\right)^\top\bfLambda_{k+1,i} + \\
&+& \left\{\sum_{i=1}^s\left[\dfrac{\partialsp\bfn_{k,i}}{\partialsp\bfy_{k}}^\top + \left(\dfrac{\partialsp e^{c_i\tau_k\bfL_{k}(\bfy_{k})}\bfy_{k}^{\mathrm{fixed}}}{\partialsp\bfy_{k}}^\top + \right.\right.\right. \nonumber\\
&+& \left.\left. \tau_k\sum\limits^{i-1}_{j=1}\dfrac{\partial(a_{ij}(\tau_k\bfL_{k}(\bfy_{k}))\bfY_{k+1,j})}{\partialsp\bfy_{k}}^\top\right)\left(\dfrac{\partialsp\bfn_{k,i}}{\partialsp\bfy}\right)^\top\right]\bfLambda_{k+1,i} + \nonumber\\
&+& \left. \left({\dfrac{\partialsp (e^{\tau_k\bfL_{k}(\bfy_{k})}\bfy_{k}^{\mathrm{fixed}})}{\partialsp\bfy_{k}}}^\top + \tau_k\sum_{i=1}^s{\dfrac{\partial(b_{i}(\tau_k\bfL_{k}(\bfy_{k}))\bfY_{k+1,i})}{\partialsp\bfy_{k}}}^\top\right)\bflambda_{k+1}\right\} \nonumber
\label{eqn:ETDRK_adjoint_update}
\end{eqnarray}
where the terms in braces are ignored if $\bfL_k$ does not depend on $\bfy_k$.
\item For $i = s, s-1, \ldots,1$, compute 
\begin{eqnarray}
\bfLambda_{k,i} &=& \bfTheta_{k,i} + 
\label{eqn:ETDRK_adjoint_internal_stages} \\
&+& \tau_{k-1}\left(\sum_{j=i+1}^s\hs a_{ji}\tau_{k-1}\bfL_{k-1}^\top\left(\dfrac{\partialsp\bfn_{k-1,j}}{\partialsp\bfy}\right)^\top\bfLambda_{k,j} + b_{i}\tau_{k-1}\bfL_{k-1}^\top\bflambda_{k}\right).\nonumber
\end{eqnarray}
The derivatives in the braces in \eqref{eqn:ETDRK_adjoint_update} are computed as discussed in Section \ref{sec:derivatives_of_phi}.
\label{8.5}
\end{subequations}
\end{itemize} \qed 
\end{algorithm}


Incidentally, inspecting the different ways in which the $\varphi_{\ell}$-functions are evaluated reveals that $\varphi_{\ell}(\tau_k\bfL_{k})^\top = \varphi_{\ell}(\tau_k\bfL_{k}^\top)$, and since in many applications the linear operator is symmetric we will have $\varphi_{\ell}(\tau_k\bfL_{k})^\top = \varphi_{\ell}(\tau_k\bfL_{k})$.


\subsection{Application to Krogstad's scheme}

We continue the example from Section~\ref{sec7.2}
and apply the above algorithm to 
a specific scheme. Recall the matrix coefficients $a_{ij}$ and $b_i$ given in \eqref{eqn:ETDRK_Butcher_tableau_Krogstad}. 
The adjoint scheme is given in Algorithm \ref{alg:adjoint_Krogstad_scheme}. 

\begin{algorithm}
\label{alg:adjoint_Krogstad_scheme}
{\rm The Adjoint Krogstad Scheme} \\[0.5em]
Let $\widehat{\bfLambda}_{k,i} = \left(\dfrac{\partialsp\bfn_{k-1,i}}{\partialsp\bfy}\right)^\top\bfLambda_{k,i}$. 
For $k = K, \ldots,1$:
\begin{itemize}
\item In this step, let $\varphi_{\ell} = \varphi_{\ell}(\tau_k\bfL_{k}(\bfy_{k}))$ and 
$\varphi_{\ell,i} = \varphi_{\ell}(c_{i}\tau_k\bfL_{k}(\bfy_{k}))$. Compute
\begin{align*}
\bflambda_{k} &= \bftheta_{k} + e^{\tau_k\bfL_{k}^\top}\left(\bflambda_{k+1} + \widehat{\bfLambda}_{k+1,4}\right) + 
e^{\frac{1}{2}\tau_k\bfL_{k}^\top}\left(\widehat{\bfLambda}_{k+1,2} + 
\widehat{\bfLambda}_{k+1,3}\right) + \widehat{\bfLambda}_{k+1,1} + \\
 & \quad + \left\{\sum\limits_{i=1}^4 \dfrac{\partialsp\bfn_{k,i}}{\partialsp\bfy_{k}}^\top\bfLambda_{k+1,i}\right. + \\
& \quad + \dfrac{\partialsp e^{\frac{1}{2}\tau_k\bfL_{k}(\bfy_{k})}\bfy_{k}^{\mathrm{fixed}}}{\partialsp\bfy_{k}}^\top\left(\widehat{\bfLambda}_{k+1,2} + \widehat{\bfLambda}_{k+1,3}\right) + \dfrac{\partialsp e^{\tau_k\bfL_{k}(\bfy_{k})}\bfy_{k}^{\mathrm{fixed}}}{\partialsp\bfy_{k}}^\top\widehat{\bfLambda}_{k+1,4} + 
 \\
& \qquad + {\dfrac{\partialsp (e^{\tau_k\bfL_{k}(\bfy_{k})}\bfy_{k}^{\mathrm{fixed}})}{\partialsp\bfy_{k}}}^\top + \frac{\tau_k}{2}\dfrac{\partial(\varphi_{1,2}\bfY_{k+1,1})}{\partialsp\bfy_{k}}^\top\widehat{\bfLambda}_{k+1,2} +  \\
& \quad + \tau_k\left(\left(\frac{1}{2}\dfrac{\partial(\varphi_{1,3}\bfY_{k+1,1})}{\partialsp\bfy_{k}}^\top + \dfrac{\partial(\varphi_{2,3}(\bfY_{k+1,2}-\bfY_{k+1,1}))}{\partialsp\bfy_{k}}^\top\right)\widehat{\bfLambda}_{k+1,3} + \right. \\
& \qquad + \left(\dfrac{\partial(\varphi_{1,4}\bfY_{k+1,1})}{\partialsp\bfy_{k}}^\top + 2\dfrac{\partial(\varphi_{2,4}(\bfY_{k+1,3}-\bfY_{k+1,1}))}{\partialsp\bfy_{k}}^\top\right)\widehat{\bfLambda}_{k+1,4} + \\
& \qquad + \left({\dfrac{\partial(\varphi_{1}\bfY_{k+1,1})}{\partialsp\bfy_{k}}}^\top - {\dfrac{\partial(\varphi_{2}\left(3\bfY_{k+1,1} - 2\bfY_{k+1,2} - 2\bfY_{k+1,3} + \bfY_{k+1,4}\right))}{\partialsp\bfy_{k}}}^\top + \right. \\
& \qquad\qquad\qquad\qquad + \left.\left.\left.4 {\dfrac{\partial(\varphi_{3}(\bfY_{k+1,1} - \bfY_{k+1,2} - \bfY_{k+1,3} + \bfY_{k+1,4}))}{\partialsp\bfy_{k}}}^\top\right)\bflambda_{k+1}\right)\right\}
\end{align*}
with $\bflambda_{K+1} = \widehat{\bfLambda}_{K+1,i} = \zero_{N \times 1}$. Ignore the terms in braces if $\bfL_k$
is independent of $\bfy_k$.
\item Now let $\varphi_{\ell} = \varphi_{\ell}(\tau_{k-1}\bfL_{k-1}(\bfy_{k-1}))$ and $\varphi_{\ell,i} = 
\varphi_{\ell}(c_{i}\tau_{k-1}\bfL_{k-1}(\bfy_{k-1}))$. The internal stages are computed as follows:
\begin{equation}
\begin{aligned}
\bfLambda_{k,4} &= \tau_{k-1}\hs b_{4}^\top\bflambda_{k} = \tau_{k-1}\left(4\varphi_{3}^\top - \varphi_{2}^\top\right)\bflambda_{k} \\
\bfLambda_{k,3} &= \tau_{k-1}\hs b_{3}^\top\bflambda_{k} + \tau_{k-1}\hs a_{4,3}^\top\widehat{\bfLambda}_{k,4} = \tau_{k-1}\hs \left(2\varphi_{2}^\top - 4\varphi_{3}^\top\right)\bflambda_{k} + \tau_{k-1}\hs 2\varphi_{2,4}^\top\widehat{\bfLambda}_{k,4} \\
\bfLambda_{k,2} &= \tau_{k-1}\hs b_{2}^\top\bflambda_{k} + \tau_{k-1}\hs a_{4,2}^\top\widehat{\bfLambda}_{k,4} + \tau_{k-1}\hs a_{3,2}^\top\widehat{\bfLambda}_{k,3} \\
&= \tau_{k-1}\left(2\varphi_{2}^\top - 4\varphi_{3}^\top\right)\bflambda_{k} + \tau_{k-1}\varphi_{2,3}^\top\widehat{\bfLambda}_{k,3} \\
\bfLambda_{k,1} &= \tau_{k-1}\hs b_{1}^\top\bflambda_{k} + \tau_{k-1}\hs a_{4,1}^\top \widehat{\bfLambda}_{k,4} + \tau_{k-1}\hs a_{3,1}^\top \widehat{\bfLambda}_{k,3} + \tau_{k-1}\hs a_{2,1}^\top \widehat{\bfLambda}_{k,2} \\
 &= \tau_{k-1}\left(\varphi_{1}^\top - 3\varphi_{2}^\top + 4\varphi_{3}^\top\right)\bflambda_{k} + \tau_{k-1}\left(\varphi_{1,4}^\top - 2\varphi_{2,4}^\top\right)\widehat{\bfLambda}_{k,4} + \\
 & \qquad\qquad\qquad\qquad\qquad + \frac{\tau_{k-1}}{2}\left(\varphi_{1,3}^\top - 2\varphi_{2,3}^\top\right) \widehat{\bfLambda}_{k,3} + \frac{\tau_{k-1}}{2}\varphi_{1,2}^\top \widehat{\bfLambda}_{k,2}.
\end{aligned}
\end{equation}
\end{itemize} \qed \\[0.7em]
\end{algorithm}
Other schemes can be handled similarly. We note that the procedure gives ample opportunity for parallelization and precomputing quantities, for instance when computing the products of $\bflambda_k$ with $\varphi_{\ell}^\top$ or the products of $\widehat{\bfLambda}_{k,i}$ with $\varphi_{\ell,i}^\top$. 



\section{Numerical Example}
\label{sec:numerical_example}

As an interesting and simple application of the methods derived above we consider the parameter estimation problem applied to the following version of the fourth-order Swift-Hohenberg model on the torus $\Omega = \mathbb{T}^2: [0,L_x) \times [0,L_y)$:
\begin{equation}
	\dfrac{\partialsp y}{\partialsp t} = r\hs y - \left(1 + \nabla^2\right)^2 y + g\hs y^2 - y^3,
\end{equation}
where $r > 0$
and $g$
are parameter functions that determine the behaviour of the solution. After spatially discretizing in some appropriate way we have
\begin{equation}
	\dfrac{\partialsp\bfy}{\partialsp t} = \mdiag{\bfr}\bfy - \left(\bfI + \nabla^2_h\right)^2 \bfy + \mdiag{\bfg}\bfy^2 - \bfy^3,
\end{equation}
with $\bfy = \bfy(t)$, $\bfr$ and $\bfg$ being the spatial discretizations of $r$ and $g$ respectively, and $\nabla^2_h$ the spatially discretized Laplace operator. The Swift-Hohenberg model is an example of a PDE whose solutions exhibit pattern formation, a phenomenon that occurs in many different branches of science, for instance biology (morphogenesis, vegetation patterns, animal markings, growth of bacterial colonies, etc.), physics (liquid crystals, nonlinear waves, B{\'e}nard cells, etc.) and chemical kinetics (e.g. the Belousov-Zhabotinsky, CIMA and PA-MBO relations). 
The field is enormous; see \cite{GollubLanger1999, MainiPainterChau1997} for some interesting applications, but there are many others. The Swift-Hohenberg model itself was derived from the equations of thermal convection \cite{SwiftHohenberg1977}, and it is possible to use other nonlinear terms than the one used here. 

The fourth-order derivative term implies that this equation is very stiff, making it a good candidate for use with an exponential integrator. The obtained patterns depend on the parameters 
$\bfm = (\bfr , \bfg )$ reshaped as a vector.
It is possible to obtain different patterns in different regions of the domain if these parameters exhibit spatial variability, 
which is the case that we consider here. This is therefore a distributed parameter estimation problem.

We take the linear part to be $\bfL = -\left(\bfI + \nabla^2_h\right)^2$,
hence $\bfn(\bfy,\bfm,t) = \mdiag{\bfr}\bfy(t) + \mdiag{\bfg}\bfy(t)^2 - \bfy(t)^3$. Notice that we have included the linear term $\mdiag{\bfr}\bfy$ in the nonlinear part of the equation: the dominant differential term is in $\bfL$ anyway.

\subsection{Experiment setup}

Using a finite difference spatial discretization with $N_x \times N_y$ grid points, the periodic boundary conditions allow us to diagonalize the linear term using the pseudospectral method, where we compute the product with the linear part in the frequency domain and then switch back to the real domain to compute the nonlinear part. Hence
\begin{equation}
	\dfrac{\partialsp\widehat{\bfy}(t)}{\partialsp t} = \widehat{\bfL}\hs\widehat{\bfy}(t) + \bfF\left(\mdiag{\bfr}\bfF^{-1}\widehat{\bfy}(t) + \mdiag{\bfg}\left(\bfF^{-1}\widehat{\bfy}(t)\right)^2 - \left(\bfF^{-1}\widehat{\bfy}(t)\right)^3\right),	\label{eqn:Swift_Hohenberg_pseudospectral}
\end{equation}
where $\bfF$ represents the 2D Fourier transform in space, $\bfF^{-1}$ is its inverse, $\widehat{\bfy}(t) = \bfF\hs\bfy(t)$ and $\widehat{\bfL}$ is the diagonalized differential operator. The wave numbers on this grid are $\bfk_{\bfx} = \frac{2\pi}{L_x}(-\frac{N_x}{2}:\frac{N_x}{2}-1)$ and $\bfk_{\bfy} = \frac{2\pi}{L_y}(-\frac{N_y}{2}:\frac{N_y}{2}-1)$, so $\nabla^2_h$ is then diagonalized with each diagonal entry the negative of the sum of the square of an element of $\bfk_{\bfx}$ and the square of an element of $\bfk_{\bfy}$. Now $\widehat{\bfy}(t)$ is taken to be the forward solution instead of $\bfy(t)$.

\begin{figure}
\centering
\begin{subfigure}{.48\textwidth}
\centering
\includegraphics[scale=0.62]{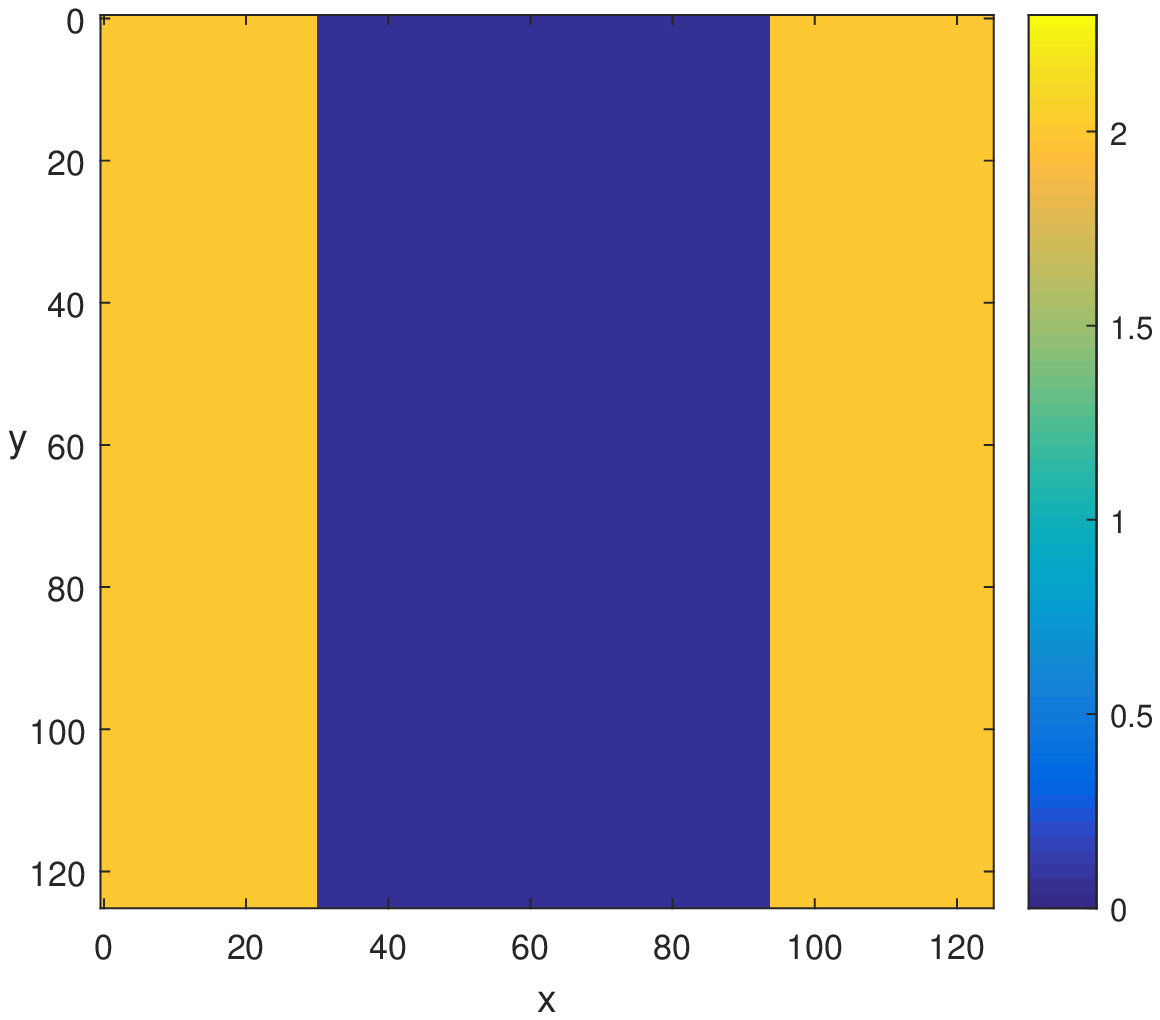}
\caption{$\bfr$}
\label{fig:actual_r}
\end{subfigure}
\begin{subfigure}{.48\textwidth}
\centering
\includegraphics[scale=0.62]{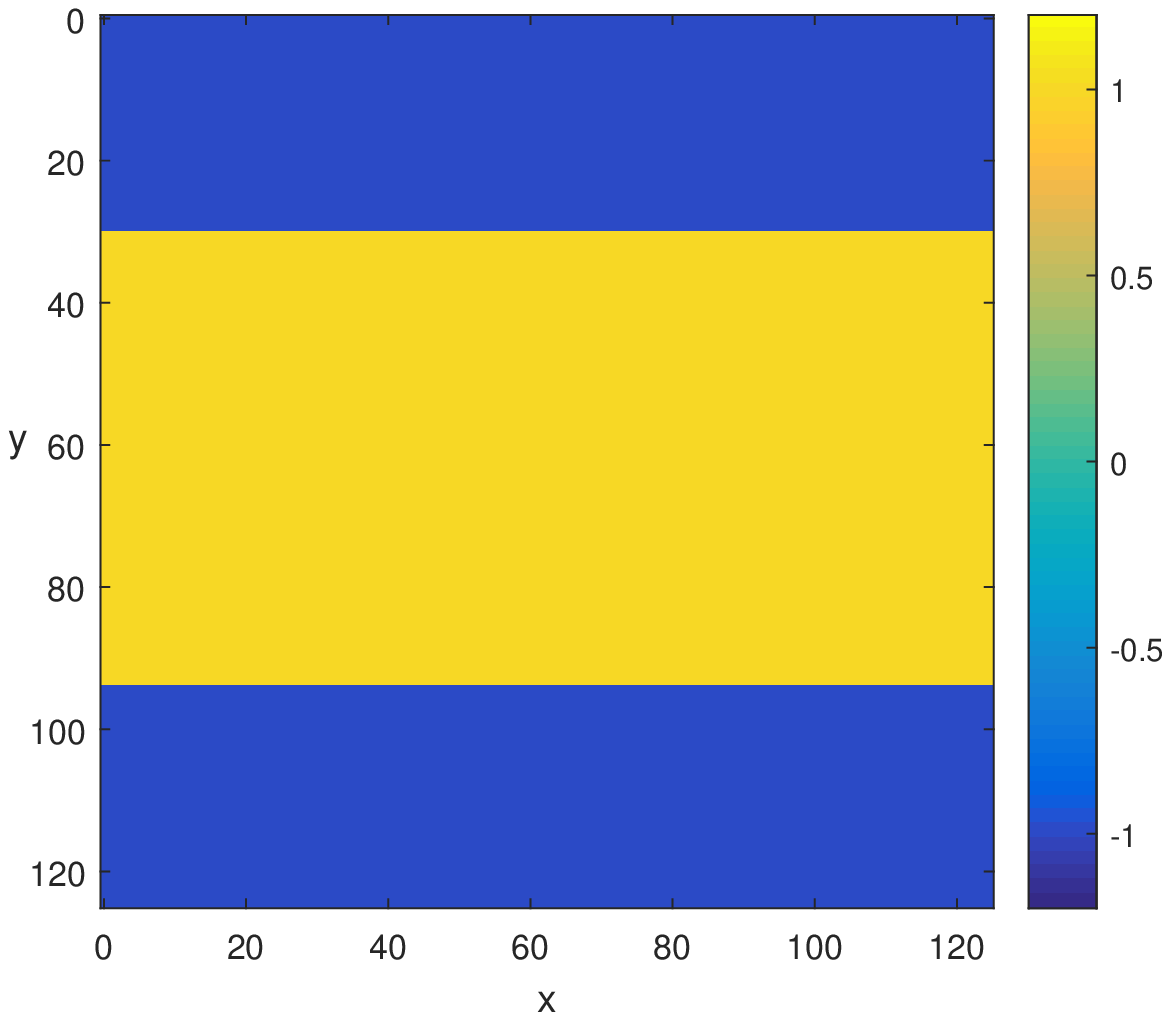}
\caption{$\bfg$}
\label{fig:actual_g}
\end{subfigure}
\caption{Actual parameter values}
\label{fig:actual_parameter_values}
\end{figure}

\begin{figure}
\includegraphics[scale=0.7]{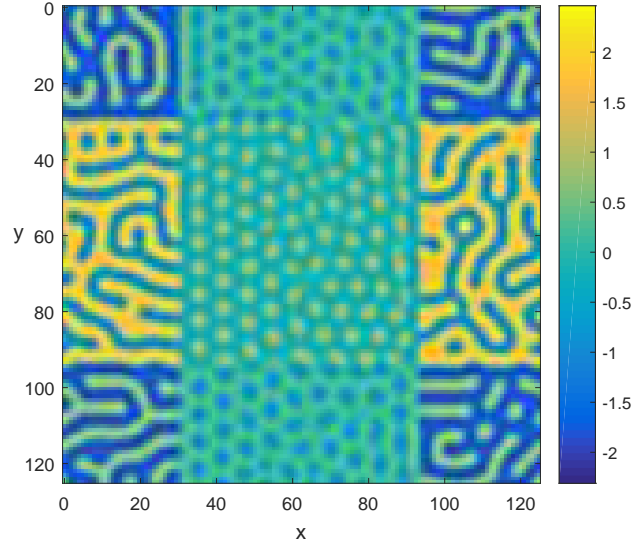}
\caption{Ground truth solution at $t = 50s$.}
\label{fig:y_at_50s}
\end{figure}

In our experiment we let $L_x = L_y = 40\pi$, $N_x = N_y = 2^7$, and the initial condition $\bfy_0$ is given by a field of Gaussian noise. We employ contour integration to evaluate the $\varphi$-functions, using a parabolic contour with 32 quadrature points. The actual parameter fields $\bfr$ and $\bfg$ are taken to be piecewise constant, as shown in Figure~\ref{fig:actual_parameter_values}. 
The value of $\bfr$ is $2$ in the outer strips and $0.04$ in the inner strip.
The value of $\bfg$ is $-1$ in the outer strips and $1$ in the inner strip.
The solution $\bfy (t)$ at time $t=50s$ is shown in Figure \ref{fig:y_at_50s}.
The outer vertical strips in the solution evolve fairly quickly relative to the central vertical strip, which evolves on a much slower time scale due to the small value of $\bfr$ there.


\subsection{Derivatives of $\bfn$}

The adjoint solution procedure and sensitivity computations require the derivatives of
$$
\bfn(\widehat{\bfy},\bfm,t) = \bfF\left(\mdiag{\bfr}\bfF^{-1}\widehat{\bfy}(t) + \mdiag{\bfg}\left(\bfF^{-1}\widehat{\bfy}(t)\right)^2 - \left(\bfF^{-1}\widehat{\bfy}(t)\right)^3\right)
$$
with respect to $\widehat{\bfy}(t)$ and $\bfm$, as well as their transposes. We have
\begin{align*}
	\dfrac{\partialsp\bfn(t)}{\partialsp\widehat{\bfy}} &= \bfF\left(\mdiag{\bfr} + 2\hs\mdiag{\bfg}\hs\left(\bfF^{-1}\widehat{\bfy}(t)\right) - 3\left(\bfF^{-1}\widehat{\bfy}(t)\right)^2\right)\bfF^{-1} \\
\dfrac{\partialsp\bfn(t)}{\partialsp\bfm} &= \begin{bmatrix} \mdiag{\widehat{\bfy}(t)} & \bfF\hs\mdiag{\left(\bfF^{-1}\widehat{\bfy}(t)\right)^2} \end{bmatrix} \\
	\dfrac{\partialsp\bfn(t)}{\partialsp\widehat{\bfy}}^\top &= \bfF^{-\top}\left(\mdiag{\bfr} + 2\hs\mdiag{\bfg}\hs\left(\bfF^{-1}\widehat{\bfy}(t)\right) - 3\left(\bfF^{-1}\widehat{\bfy}(t)\right)^2\right)\bfF^{\top} \\
\dfrac{\partialsp\bfn(t)}{\partialsp\bfm}^\top &= \begin{bmatrix} \mdiag{\widehat{\bfy}(t)} \\ \mdiag{\left(\bfF^{-1}\widehat{\bfy}(t)\right)^2}\bfF^\top \end{bmatrix}.
\end{align*}
The transposes of $\bfF$ and $\bfF^{-1}$ are $\bfF^\top = N\bfF^{-1}$ and $\bfF^{-\top} = \frac{1}{N}\bfF$, with $N = N_xN_y$, so that
\begin{align*}
	\dfrac{\partialsp\bfn(t)}{\partialsp\widehat{\bfy}}^\top &= \bfF\left(\mdiag{\bfr} + 2\hs\mdiag{\bfg}\hs\left(\bfF^{-1}\widehat{\bfy}(t)\right) - 3\left(\bfF^{-1}\widehat{\bfy}(t)\right)^2\right)\bfF^{-1} \\
\dfrac{\partialsp\bfn(t)}{\partialsp\bfm}^\top &= \begin{bmatrix} \mdiag{\widehat{\bfy}(t)} \\ N\mdiag{\left(\bfF^{-1}\widehat{\bfy}(t)\right)^2}\bfF^{-1} \end{bmatrix}.
\end{align*}
Note that $\dfrac{\partialsp\bfn(t)}{\partialsp\widehat{\bfy}}^\top = \dfrac{\partialsp\bfn(t)}{\partialsp\widehat{\bfy}}$.


\subsection{Order of Accuracy}

We numerically illustrate that the adjoint solution and the gradient have the same order of accuracy as the forward solution, meaning that $p$th-order forward schemes are expected to lead to $p$th-order adjoint schemes. In this subsection we let $\bfy(n\tau)$, $\bflambda(n\tau)$ and $\nabla(n\tau)$ denote the quantity of interest one gets when performing the computations with a time-step $n\tau$, where we let $\tau = \frac{1}{80}s$. We let $\bfy_{\mathrm{exact}}$, $\bflambda_{\mathrm{exact}}$ and $\nabla_{\mathrm{exact}}$ denote the "exact" values attained by performing the computations with a fine time-step $\frac{1}{160}s$ using the Krogstad scheme. The error between the computed and exact forward solution is $\epsilon_{\bfy}(n\tau) = \bfy(n\tau) - \bfy_{\mathrm{exact}}$, and we must have $\|\epsilon_{\bfy}(n\tau)\|\approx O(n\tau)^p$ for a $p$th-order method, from which it follows that $\left(\|\epsilon_{\bfy}(2^{i}\tau)\|/\|\epsilon_{\bfy}(2^{i-1}\tau)\|\right) \approx 2^p$. Computing 
$p_{\bfy}\left(2^{i}\tau,2^{i-1}\tau\right) := \log_2\left(\|\epsilon_{\bfy}(2^{i}\tau)\|/\|\epsilon_{\bfy}(2^{i-1}\tau)\|\right)$ for different values of $i = 1,2,\ldots$ then gives an approximation of the order of accuracy $p$, with the estimate being more accurate for smaller values of $\tau$ and $i$. 
\begin{table}
\small
\centering
\begin{tabular}{c | C{6em} C{6em} C{6em} C{6em}}
	& Euler & Cox-Matthews & Krogstad & Hochbruck-Ostermann \\ \hline
	$p_{\bfy}\left(2\tau,\tau\right)$ & $0.9914$ & $4.0644$ & $4.0699$ & $4.0275$ \\
	$p_{\bfy}\left(4\tau,2\tau\right)$ & $0.9908$ & $3.9726$ & $3.9732$ & $3.9719$ \\
	$p_{\bfy}\left(8\tau,4\tau\right)$ & $-$ & $3.9375$ & $3.9378$ & $3.9387$ \\
	$p_{\bfy}\left(16\tau,8\tau\right)$ & $-$ & $3.8849$ & $3.8835$ & $3.8770$ \\ \hline
\end{tabular}
\caption{Approximate order of accuracy $p_{\bfy}$ for the forward solution, computed using $p_{\bfy}\left(2^{i}\tau,2^{i-1}\tau\right) = \log_2\left(\|\epsilon_{\bfy}(2^{i}\tau)\|/\|\epsilon_{\bfy}(2^{i-1}\tau)\|\right)$ for $i = 1,2,3,4$, with $\epsilon_{\bfy}(2^{i}\tau) = \bfy(2^{i}\tau) - \bfy_{\mathrm{exact}}$.}
\label{tab:order_of_accuracy_y}
\end{table}
The results for the four ETD schemes mentioned in this paper are given in Table \ref{tab:order_of_accuracy_y}. 
We note the following:
\begin{itemize}
\item The simulations were run from 0s to 20s for the forward solution, and from 20s to 0s for the adjoint solution.
\item The initial condition is Gaussian noise and the approximate orders of accuracy actually differ from simulation to simulation because of this. These differences are only slight for the higher-order methods, but can be quite pronounced for the lower-order Euler method. Therefore the results shown are the averages of $10$ simulations using different initial conditions.
\item The larger time steps are too large for the Euler method, and even the smaller time-steps can lead to inaccurate results on occasion. In computing the average order of accuracy we have therefore only included the computed $p_{\bfy}$'s that are in the interval $(p-0.5,p+0.5)$.
\item Incidentally, for the Swift-Hohenberg equation with periodic boundary conditions, we see that Cox-Matthews and Krogstad schemes do indeed attain fourth order accuracy.
\end{itemize}
The quantities $p_{\bflambda}\left(2^{i}\tau,2^{i-1}\tau\right)$ and $p_{\nabla}\left(2^{i}\tau,2^{i-1}\tau\right)$ are defined analogously to $p_{\bfy}$ and are given in Tables~\ref{tab:order_of_accuracy_lambda} and~\ref{tab:order_of_accuracy_grad}, respectively. We have again averaged the approximate orders of accuracy from 10 simulations.
\begin{table}
\small
\centering
\begin{tabular}{c | C{6em} C{6em} C{6em} C{6em}}
	& Euler & Cox-Matthews & Krogstad & Hochbruck-Ostermann \\ \hline
	$p_{\bflambda}\left(2\tau,\tau\right)$ & $0.9976$ & $4.0383$ & $4.0588$ & $3.9902$ \\
	$p_{\bflambda}\left(4\tau,2\tau\right)$ & $0.9434$ & $3.9516$ & $3.9568$ & $3.9679$ \\
	$p_{\bflambda}\left(8\tau,4\tau\right)$ & $-$ & $3.8969$ & $3.9041$ & $3.9343$ \\
	$p_{\bflambda}\left(16\tau,8\tau\right)$ & $-$ & $3.8027$ & $3.8100$ & $3.8616$ \\ \hline
\end{tabular}
\caption{Approximate order of accuracy $p_{\bflambda}$ for the adjoint solution, computed using $p_{\bflambda}\left(2^{i}\tau,2^{i-1}\tau\right) = \log_2\left(\|\epsilon_{\bflambda}(2^{i}\tau)\|/\|\epsilon_{\bflambda}(2^{i-1}\tau)\|\right)$ for $i = 1,2,3,4$, with $\epsilon_{\bflambda}(2^{i}\tau) = \bflambda(2^{i}\tau) - \bflambda_{\mathrm{exact}}$.}
\label{tab:order_of_accuracy_lambda}
\end{table}

\begin{table}
\small
\centering
\begin{tabular}{c | C{6em} C{6em} C{6em} C{6em}}
	& Euler & Cox-Matthews & Krogstad & Hochbruck-Ostermann \\ \hline
	$p_{\nabla}\left(2\tau,\tau\right)$ & $1.0260$ & $4.0430$ & $4.0627$ & $3.9947$ \\
	$p_{\nabla}\left(4\tau,2\tau\right)$ & $0.9270$ & $3.9546$ & $3.9575$ & $3.9684$ \\
	$p_{\nabla}\left(8\tau,4\tau\right)$ & $-$ & $3.9022$ & $3.9056$ & $3.9347$ \\
	$p_{\nabla}\left(16\tau,8\tau\right)$ & $-$ & $3.8103$ & $3.8136$ & $3.8622$ \\ \hline
\end{tabular}
\caption{Approximate order of accuracy $p_{\nabla}$ for the gradient, computed using $p_{\nabla}\left(2^{i}\tau,2^{i-1}\tau\right) = \log_2\left(\|\epsilon_{\nabla}(2^{i}\tau)\|/\|\epsilon_{\nabla}(2^{i-1}\tau)\|\right)$ for $i = 1,2,3,4$, with $\epsilon_{\nabla}(2^{i}\tau) = \nabla(2^{i}\tau) - \nabla_{\mathrm{exact}}$.}
\label{tab:order_of_accuracy_grad}
\end{table}
The crucial observation here is that the orders of accuracy of the adjoint ETD schemes and the resulting gradient are the same as that of the corresponding forward scheme.

\subsection{Testing the Rosenbrock Approach}

To test the order of accuracy of adjoint exponential Rosenbrock methods, and simulataneously check that our derivations and implementations for these methods are correct, the Swift-Hohenberg equation was reformulated by finding the Jacobian of the right-hand side of \eqref{eqn:Swift_Hohenberg_pseudospectral},
\[
	\dfrac{\partialsp\bff}{\partialsp\widehat{\bfy}} = \widehat{\bfL} + \bfF\left(\mdiag{\bfr} + 2\hs\mdiag{\bfg}\mdiag{\bfF^{-1}\widehat{\bfy}(t)} - 3\hs\mdiag{\bfF^{-1}\widehat{\bfy}(t)}\right)\bfF^{-1},
\]
and then setting $\bfL_k = \dfrac{\partialsp\bff(\widehat{\bfy}_k)}{\partialsp\widehat{\bfy}}$ at the $k$th time-step. Consequently, with $\bfy_k = \bfF^{-1}\widehat{\bfy}_k$,
\[
	\bfL_k = \widehat{\bfL} + \bfF\mdiag{\bfr + 2\hs\bfg\odot\bfy_k - 3\hs\bfy_k^2}\bfF^{-1},
\]
and hence
\[
	\bfn_k = \bff - \bfL_k\hs\widehat{\bfy}(t) = \bfF\hs\mdiag{\bfg\odot\left(\bfF^{-1}\widehat{\bfy}(t) - 2\bfy_k\right) - \left(\bfF^{-1}\widehat{\bfy}(t)\right)^2 - 3\hs\bfy_k^2}\bfF^{-1}\widehat{\bfy}(t).
\]
\begin{table}
\small
\centering
\begin{tabular}{c | C{6em} C{6em} C{6em} C{6em}}
	& Euler & Cox-Matthews & Krogstad & Hochbruck-Ostermann \\ \hline
	$p_{\bfy}\left(4\tau,2\tau\right)$ & $1.9505$ & $4.1970$ & $4.1967$ & $4.1213$ \\
	$p_{\bfy}\left(8\tau,4\tau\right)$ & $1.9713$ & $4.0774$ & $4.0811$ & $4.0561$ \\
	$p_{\bfy}\left(16\tau,8\tau\right)$ & $1.9846$ & $3.7946$ & $3.7991$ & $3.7884$ \\ \hline
\end{tabular}
\caption{Approximate order of accuracy $p_{\bfy}$ for the forward solution using a Rosenbrock approach, computed using $p_{\bfy}\left(2^{i}\tau,2^{i-1}\tau\right) = \log_2\left(\|\epsilon_{\bfy}(2^{i}\tau)\|/\|\epsilon_{\bfy}(2^{i-1}\tau)\|\right)$ for $i = 1,2,3,4$, with $\epsilon_{\bfy}(2^{i}\tau) = \bfy(2^{i}\tau) - \bfy_{\mathrm{exact}}$.}
\label{tab:order_of_accuracy_y_rosenbrock}
\end{table}
\begin{table}
\small
\centering
\begin{tabular}{c | C{6em} C{6em} C{6em} C{6em}}
	& Euler & Cox-Matthews & Krogstad & Hochbruck-Ostermann \\ \hline
	$p_{\bflambda}\left(4\tau,2\tau\right)$ & $1.8885$ & $4.2572$ & $4.2482$ & $4.2244$ \\
	$p_{\bflambda}\left(8\tau,4\tau\right)$ & $1.9414$ & $4.0531$ & $4.0498$ & $4.0411$ \\
	$p_{\bflambda}\left(16\tau,8\tau\right)$ & $1.9702$ & $3.1043$ & $3.1029$ & $3.0998$ \\ \hline
\end{tabular}
\caption{Approximate order of accuracy $p_{\bflambda}$ for the adjoint solution using a Rosenbrock approach, computed using $p_{\bflambda}\left(2^{i}\tau,2^{i-1}\tau\right) = \log_2\left(\|\epsilon_{\bflambda}(2^{i}\tau)\|/\|\epsilon_{\bflambda}(2^{i-1}\tau)\|\right)$ for $i = 1,2,3,4$, with $\epsilon_{\bflambda}(2^{i}\tau) = \bflambda(2^{i}\tau) - \bflambda_{\mathrm{exact}}$.}
\label{tab:order_of_accuracy_lambda_rosenbrock}
\end{table}
\begin{table}
\small
\centering
\begin{tabular}{c | C{6em} C{6em} C{6em} C{6em}}
	& Euler & Cox-Matthews & Krogstad & Hochbruck-Ostermann \\ \hline
	$p_{\bflambda}\left(4\tau,2\tau\right)$ & $1.8730$ & $4.2892$ & $4.2779$ & $4.2580$ \\
	$p_{\bflambda}\left(8\tau,4\tau\right)$ & $1.9345$ & $4.0911$ & $4.0839$ & $4.0764$ \\
	$p_{\bflambda}\left(16\tau,8\tau\right)$ & $1.9669$ & $3.4815$ & $3.4758$ & $3.4754$ \\ \hline
\end{tabular}
\caption{Approximate order of accuracy $p_{\nabla}$ for the gradient using a Rosenbrock approach, computed using $p_{\nabla}\left(2^{i}\tau,2^{i-1}\tau\right) = \log_2\left(\|\epsilon_{\nabla}(2^{i}\tau)\|/\|\epsilon_{\nabla}(2^{i-1}\tau)\|\right)$ for $i = 1,2,3,4$, with $\epsilon_{\nabla}(2^{i}\tau) = \nabla(2^{i}\tau) - \nabla_{\mathrm{exact}}$.}
\label{tab:order_of_accuracy_grad_rosenbrock}
\end{table}
The experiment from the previous subsection is repeated, but due to the significant increase in computational effort required by the additional terms in the derivatives we have run the simulations for only 2 different random initial conditions, and for a smallest time-step of $2\tau$, with the "exact" equations computed using a time-step of $\tau$. We have also run the simulations for just 10s instead of 20s. The results in tables \ref{tab:order_of_accuracy_y_rosenbrock}-\ref{tab:order_of_accuracy_grad_rosenbrock} suggest that the adjoint exponential Rosenbrock method does indeed also attain the same (numerical) order of accuracy as the corresponding forward method. Oddly the Euler method seems to have an order of accuracy of around 2 instead of the expected value of 1, but this result should be viewed with reservation. The results for the largest time-step suggest that this time-step was too large, with a noticeable decrease in the order of accuracy especially for the adjoint solution.

\subsection{Parameter Estimation}

\begin{figure}
\centering
\begin{subfigure}{.48\textwidth}
\centering
\includegraphics[scale=0.62]{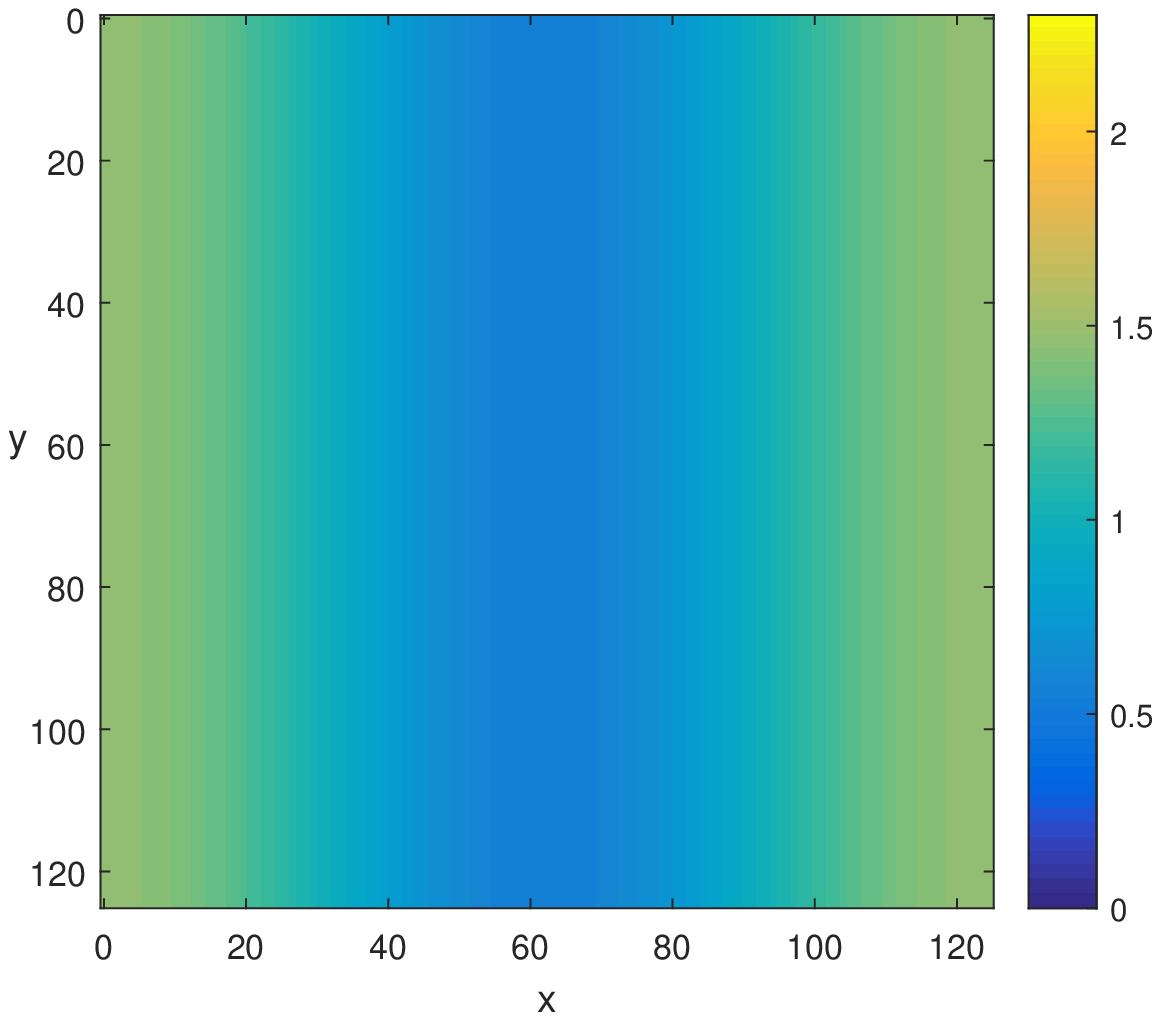}
\caption{Initial guess for $\bfr$}
\label{fig:r_initial}
\end{subfigure}
\begin{subfigure}{.48\textwidth}
\centering
\includegraphics[scale=0.62]{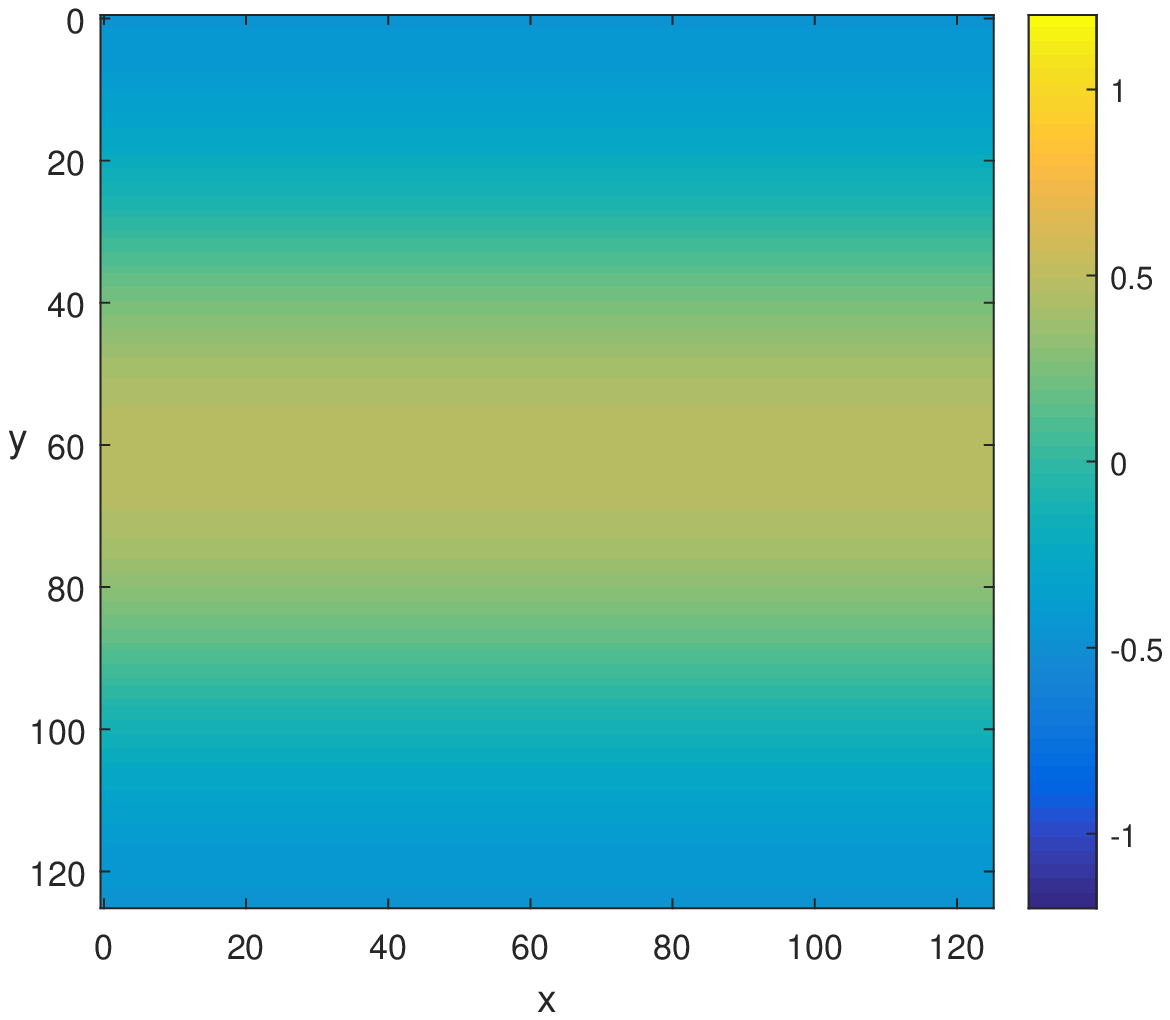}
\caption{Initial guess for $\bfg$}
\label{fig:g_initial}
\end{subfigure} \\
\begin{subfigure}{.48\textwidth}
\centering
\includegraphics[scale=0.66]{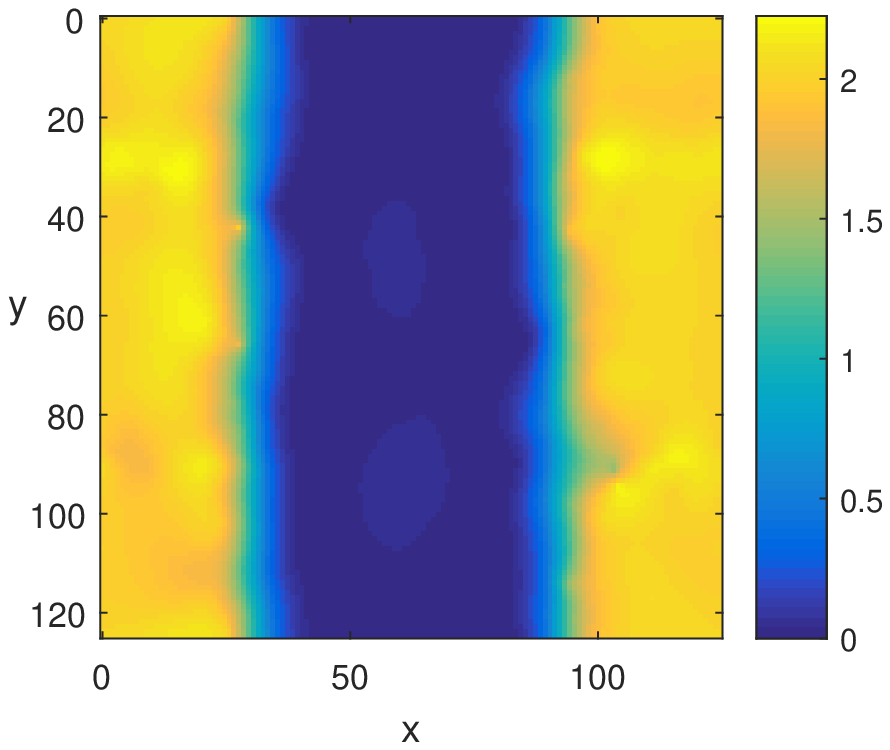}
\caption{Recovered $\bfr$}
\label{fig:r_recovered}
\end{subfigure}
\begin{subfigure}{.48\textwidth}
\centering
\includegraphics[scale=0.66]{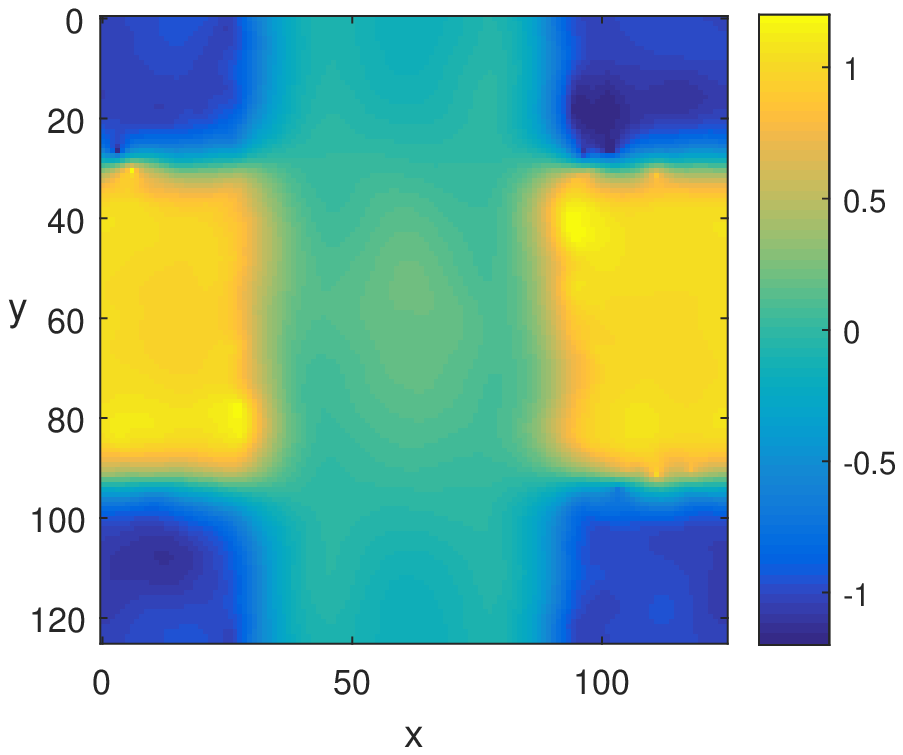}
\caption{Recovered $\bfg$}
\label{fig:g_recovered}
\end{subfigure}
\caption{Recovered parameter values (bottom) and initial guesses (top)}
\label{fig:initial_and_recovered_parameters}
\end{figure}

Although outside the stated scope of this paper, we show a possible parameter estimate attained using a gradient-based optimization method, where the gradient is computed using the results from this paper. Recovering estimates that are closer to the ground truth parameters takes a more sophisticated approach than the one used here.
In particular, the regularization term $\opR(\bfm, \bfm^\mathrm{ref})$ in~\eqref{eqn:minimize_Omega}, 
which is not the focus of this paper, may have to be altered, or a level set method may be introduced; 
this is the subject of a future investigation. The initial guesses for the parameters are shown in Figures \ref{fig:r_initial} and \ref{fig:g_initial}, and the recovered parameters are shown in Figures \ref{fig:r_recovered} and \ref{fig:g_recovered}.
Here are details of the setup:
\begin{itemize}
\item The simulations were run using the Krogstad scheme with a time-step of $0.25$s.
\item Observations were taken every $0.5$s up to $25$s, and 5\% noise was added to each observation.
\item The standard least-squares misfit function $\opM = \|\bfd(\bfm) - \bfd^{\mathrm{obs}}\|$ was employed.
\item TV regularization using the smoothed Huber norm was used with $\beta = 10$.
\item The parameters were recovered using 400 iterations of L-BFGS (keeping 20 previous updates in storage) with cubic line search. The results after 200 iterations were already very similar to those in Figure \ref{fig:initial_and_recovered_parameters}.
\item We used a projected gradient method, where the value of $\bfr$ at each point was constrained to always lie in the interval $[0.01,2.3]$ and the value of $\bfg$ at each point was in $[-1.2,1.2]$.
\end{itemize}
\noindent We make the following remarks on the results:
\begin{itemize}
\item The recovered values of $\bfr$ are quite acceptable, as are the values of $\bfg$ on the two outside vertical strips, whereas the estimate of the central vertical strip of $\bfg$ is poor.
\item This is the result of the observations having only been taken for up to 25s, a time
after the patterns on the outside vertical strips have formed but before the central pattern, which evolves on a much slower time scale, has formed.
\end{itemize}



\section{Conclusions and Future Work}
\label{sec:conclusion}

This paper considers the application of the discrete adjoint method to exponential integration methods for the purposes of 
(distributed) parameter estimation and sensitivity analysis in time-dependent PDEs. 
We have derived algorithms for the linearized forward and, more importantly, the adjoint problem, and have applied these algorithms, as an instance of these integration methods, to a specific scheme.

We have found that if the linear operator depends on the solution $\bfy_k$ at the current time step 
(as is the case with {\em some} Rosenbrock-type methods) or the model parameters, the derivatives of the $\varphi$-functions with respect to $\bfy_k$ and $\bfm$ introduce significant computational overhead. 
The extra expense introduced in this case 
could be prohibitive in some applications, and it might then be more reasonable to instead apply an IMEX method to the linearized PDE. The use of IMEX methods in the context of parameter estimation will be the subject of a future investigation.

A simple experiment reveals that the adjoint exponential integrator and the computed gradient of the misfit function have the same order of accuracy as the corresponding exponential integrator. This of course does not constitute a general proof. In case that the linear operator $\bfL$ is independent of the forward solution and the model parameters it appears that one should be able to analytically prove the conjecture that this observation offers. A general proof in the case of Rosenbrock-type methods, however, appears to be a more remote prospect given how dependent the evaluations of the derivatives of the $\varphi$-functions are on their numerical implementation.

The results of this work will be of interest in applications where exponential integration is the best choice for solving the forward solution, in particular when the PDE is stiff. One such application of interest to us is pattern formation when the model parameters exhibit spatial variability, and we have applied the techniques from this paper to a simple parameter estimation problem involving Swift-Hohenberg model. A more in-depth examination into parameter estimation of pattern formation problems will be the subject of future work.


\bibliographystyle{siam}
\bibliography{bibliography}  

\begin{thebibliography}{10}

\bibitem{Certaine1960}
{\sc J.~Certaine}, {\em The solution of ordinary differential equations with
  large time constants}, in Mathematical {M}ethods for {D}igital {C}omputers,
  A.~Ralston and H.~Wilf, eds., Wiley, New York, 1960.

\bibitem{Chavent1974}
{\sc G.~Chavent}, {\em Identification of function parameters in partial
  differential equations}, in Identification of Parameters in Distributed
  Systems: Symposium, Joint Automatic Control Conference, R.~Goodson, ed.,
  American Society of Mechanical Engineers, New York, 1974.

\bibitem{CoxMatthews2002}
{\sc S.~M. Cox and P.~C. Matthews}, {\em Exponential {T}ime {D}ifferencing for
  stiff systems}, Journal of Computational Physics, 176 (2002),
  pp.~430–--455.

\bibitem{DennisSchnabel1996}
{\sc J.~E. Dennis and R.~B. Schnabel}, {\em Numerical Methods for Unconstrained
  Optimization and Nonlinear Equations}, SIAM, 1996.

\bibitem{EiermannErnst2006}
{\sc M.~Eiermann and O.~G. Ernst}, {\em A restarted {K}rylov subspace method
  for the evaluation of matrix functions}, SIAM Journal of Numerical Analysis,
  44 (2006), pp.~2481–--2504.

\bibitem{EnglHankeNeubauer1996}
{\sc H.~W. Engl, M.~Hanke, and A.~Neubauer}, {\em Regularization of Inverse
  Problems}, Kluwer, 1996.

\bibitem{GollubLanger1999}
{\sc J.~P. Gollub and J.~S. Langer}, {\em Pattern formation in nonequilibrium
  physics}, Revisions of Modern Physics, 71 (1999), pp.~S396--S403.

\bibitem{GriewankWalther2008}
{\sc A.~Griewank and A.~Walther}, {\em Evaluating Derivatives: Principles and
  Techniques of Algorithmic Differentiation}, vol.~105 of Other Titles in
  Applied Mathematics, SIAM, 2nd~ed., 2008.

\bibitem{Gunzburger2002}
{\sc M.~D. Gunzburger}, {\em Perspectives in Flow Control and Optimization},
  SIAM, 2002.

\bibitem{Guttel2013}
{\sc S.~G{\"u}ttel}, {\em Rational krylov approximation of matrix functions:
  {N}umerical methods and optimal pole selection}, GAMM Mitteilungen, 36
  (2013), pp.~8–--31.

\bibitem{HochbruckLubich1997}
{\sc M.~Hochbruck and C.~Lubich}, {\em On krylov subspace approximations to the
  matrix exponential operator}, SIAM Journal of Numerical Analysis, 34 (1997),
  pp.~1911–--1925.

\bibitem{HochbruckOstermann2005a}
{\sc M.~Hochbruck and A.~Ostermann}, {\em Explicit exponential {R}unge–
  {K}utta methods for semilinear parabolic problems}, SIAM Journal of Numerical
  Analysis, 43 (2005), pp.~1069–--1090.

\bibitem{HochbruckOstermann2010}
\leavevmode\vrule height 2pt depth -1.6pt width 23pt, {\em Exponential
  {I}ntegrators}, Acta Numerica, 19 (2010), pp.~209–--286.

\bibitem{HochbruckOstermannSchweitzer2009}
{\sc M.~Hochbruck, A.~Ostermann, and J.~Schweitzer}, {\em Exponential
  {R}osenbrocktype methods}, SIAM Journal of Numerical Analysis, 47 (2009),
  pp.~786--–803.

\bibitem{HochbruckVanDenEshof2006a}
{\sc M.~Hochbruck and J.~van~den Eshof}, {\em Explicit integrators of
  {R}osenbrock-type}, Oberwolfach Reports, 3 (2006), pp.~1107–--1110.

\bibitem{KassamTrefethen2005}
{\sc A.-K. Kassam and L.~N. Trefethen}, {\em Fourth-order time-stepping for
  stiff {PDE}s}, SIAM Journal of Scientific Computing, 26 (2005),
  pp.~1214--1233.

\bibitem{Krogstad2005}
{\sc S.~Krogstad}, {\em Generalized integrating factor methods for stiff
  {PDE}s}, Journal of Computational Physics, 203 (2005), pp.~72–--88.

\bibitem{LuanOstermann2013}
{\sc V.~T. Luan and A.~Ostermann}, {\em Exponential {B}-series: The stiff
  case}, SIAM Journal of Numerical Analysis, 51 (2013), pp.~3431--–3445.

\bibitem{LuanOstermann2014b}
\leavevmode\vrule height 2pt depth -1.6pt width 23pt, {\em Explicit exponential
  {R}unge-{K}utta methods of high order for parabolic problems}, Journal of
  Computational and Applied Mathematics, 256 (2014), pp.~168–--179.

\bibitem{LuanOstermann2014a}
\leavevmode\vrule height 2pt depth -1.6pt width 23pt, {\em Exponential
  {R}osenbrock methods of order five — construction, analysis and numerical
  comparisons}, Journal of Computational and Applied Mathematics, 255 (2014),
  pp.~417–--431.

\bibitem{LuanOstermann2014c}
\leavevmode\vrule height 2pt depth -1.6pt width 23pt, {\em Stiff order
  conditions for exponential {R}unge-{K}utta methods of order five}, in
  Modeling, Simulation and Optimization of Complex Processes - HPSC 2012,
  H.~G.~B. et~al., ed., Springer, 2014.

\bibitem{LuanOstermann2016}
\leavevmode\vrule height 2pt depth -1.6pt width 23pt, {\em Parallel exponential
  {R}osenbrock methods}, Computers and mathematics with Applications, 71
  (2016), pp.~1137--1150.

\bibitem{MainiPainterChau1997}
{\sc P.~K. Maini, K.~J. Painter, and H.~N.~P. Chau}, {\em Spatial pattern
  formation in chemical and biological systems}, Faraday Transactions, 93
  (1997), pp.~3601--3610.

\bibitem{Neidinger2010}
{\sc R.~Neidinger}, {\em Introduction to {A}utomatic {D}ifferentiation and
  {MATLAB} {O}bject-{O}riented programming}, SIAM Review, 52 (2010),
  pp.~545--–563.

\bibitem{NocedalWright2006}
{\sc J.~Nocedal and S.~J. Wright}, {\em Numerical Optimization}, Springer,
  2~ed., 2006.

\bibitem{Plessix2006}
{\sc F.-E. Plessix}, {\em A review of the adjoint-state method for computing
  the gradient of a functional with geophysical applications}, Geophys. J.
  Int., 167 (2006), pp.~495--–503.

\bibitem{Pope1963}
{\sc D.~A. Pope}, {\em An exponential method of numerical integration of
  ordinary differential equations}, Communications of the Association for
  Computing Machinery, 6 (1963), pp.~491--–493.

\bibitem{Saad2003}
{\sc Y.~Saad}, {\em Iterative {M}ethods for {S}parse {L}inear {S}ystems}, SIAM,
  2nd~ed., 2003.

\bibitem{SchmelzerTrefethen2007}
{\sc T.~Schmelzer and L.~N. Trefethen}, {\em Evaluating matrix functions for
  exponential integrators via {C}arath{\'e}odory-{F}ej{\'e}r approximation and
  contour integral}, Electronic Transactions on Numerical Analysis, 29 (2007),
  pp.~1–--18.

\bibitem{SwiftHohenberg1977}
{\sc J.~Swift and P.~Hohenberg}, {\em Hydrodynamic fluctuations at the
  convective instability}, Phys. Rev. A., 15 (1977), pp.~319–--328.

\bibitem{Tokman2005}
{\sc M.~Tokman}, {\em Efficient integration of large stiff systems of {ODE}s
  with exponential propagation iterative ({EPI}) methods}, Journal of
  Computational Physics, 213 (2005), pp.~748–--776.

\bibitem{Tokman2011}
\leavevmode\vrule height 2pt depth -1.6pt width 23pt, {\em A new class of
  exponential propagation iterative methods of {R}unge–{K}utta type
  ({EPIRK})}, Journal of Computational Physics, 230 (2011), pp.~8762--–8778.

\bibitem{TrefethenWeidemanSchmelzer2006}
{\sc L.~N. Trefethen, J.~A.~C. Weideman, and T.~Schmelzer}, {\em Talbot
  {Q}uadratures and {R}ational {A}pproximations}, BIT, 46 (2006),
  pp.~653–--670.

\bibitem{vandenDoelAscherHaber2013}
{\sc K.~van~den Doel, U.~Ascher, and E.~Haber}, {\em The lost honour of
  $\ell_2$-based regularization}, Radon Series in Computational and Applied
  Math,  (2013).
\newblock M. Cullen, M. Freitag, S. Kindermann and R. Scheinchl (Eds).

\bibitem{Vogel2002}
{\sc C.~R. Vogel}, {\em Computational {M}ethods for {I}nverse {P}roblems},
  SIAM, Philadelphia, 2002.

\end{thebibliography}



\end{document}